\documentclass[twocolumn]{autart}    

\usepackage{graphicx}          
\usepackage{amsmath}
\usepackage{amsfonts}
\usepackage{amssymb} 

\usepackage{theoremref}

\usepackage{subfigure}

\usepackage{stfloats}

\usepackage{color}

\usepackage{tcolorbox}

\begin{document}
	
	\begin{frontmatter}
		
		\title{Dominance Regions of Pursuit-evasion Games in Non-anticipative Information Patterns\thanksref{footnoteinfo}} 
		
		\thanks[footnoteinfo]{This paper was not presented at any IFAC 
			meeting. Corresponding author Fang Deng. The work was supported in part by the National Key R\&D Program of China under Grant 2022ZD0119701, the National Natural Science Foundation of China National Science Fund for Distinguished Young Scholars 62025301, the National Natural Science Foundation of China under Grant 61933002, and the National Natural Science Foundation of China Basic Science Center Program under Grant 62088101.}
		
		\author[BIT-math,Key Lab]{Weiwen Huang}\ead{huangweiwen96@hotmail.com},    
		\author[BUCT-CIST]{Li Liang}\ead{lianglibuct@buct.edu.cn},               
		\author[Key Lab,BIT-auto,chongqing]{Ningsheng Xu}\ead{Xuningsheng1@163.com},  
		\author[Key Lab,BIT-auto,chongqing]{Fang Deng}\ead{dengfang@bit.edu.cn}
		
		\address[BIT-math]{School of Mathematics and Statistics, Beijing Institute of Technology, Beijing 100081, China}  
		\address[Key Lab]{National Key Lab of Autonomous Intelligent Unmanned Systems, Beijing 100081, China}
		\address[BUCT-CIST]{College of Information Science and Technology, Beijing University of Chemical Technology, Beijing 100029, China}                                              
		\address[BIT-auto]{School of Automation, Beijing Institute of Technology, Beijing 100081, China}             
		\address[chongqing]{Beijing Institute of Technology Chongqing Innovation Center, Chongqing 401120, China}
		
		
		\begin{abstract}                          
			The evader's dominance region is an important concept and the foundation of geometric methods for pursuit-evasion games. This article mainly reveals the relevant properties of the evader's dominance region, especially in non-anticipative information patterns. 
			We can use these properties to research pursuit-evasion games in non-anticipative information patterns. The core problem is under what condition the pursuer has a non-anticipative strategy to prevent the evader leaving its initial dominance region before being captured regardless of the evader's strategy.
			We first define the evader's dominance region by the shortest path distance, and we rigorously prove for the first time that the initial dominance region of the evader is the reachable region of the evader in the open-loop sense. Subsequently, we prove that there exists a non-anticipative strategy by which the pursuer can capture the evader before the evader leaves its initial dominance region's closure in the absence of obstacles. For cases with obstacles, we provide a counter example to illustrate that such a non-anticipative strategy does not always exist, and provide a necessary condition for the existence of such strategy. Finally, we consider a scenario with a single corner obstacle and provide a sufficient condition for the existence of such a non-anticipative strategy. At the end of this article, we discuss the application of the evader's dominance region in target defense games. This article has important reference significance for the design of non-anticipative strategies in pursuit-evasion games with obstacles.
			
		\end{abstract}
		
	\end{frontmatter}
	
	\section{Introduction}
	The pursuit-evasion game provides a general framework to mathematically formalize many practical application problems. These practical problems include control, navigation, biological behavior analysis and military confrontation. In the simplest pursuit-evasion game, there are two players: one pursuer and one evader. The pursuer intends to capture or intercept the evader, and the evader intends to avoid the capture or interception of the pursuer to achieve a certain purpose. The differential game problem can be regarded as an extension of the control problem. The two-side maximum principle and Hamilton-Jacobi-Isaacs equation are important tools to solve  pursuit-evasion games \cite{MR0210469}\cite{MR1311921}. Information patterns play an important role in differential games. The non-anticipative information pattern is widely used in much literature \cite{elliott1972existence}\cite{cardaliaguet1996differential}\cite{mitchell2005time}.  It is more practical than the open-loop information pattern. Non-anticipative strategies do not rely on future information. 
	Some commonly used strategies, such as state feedback strategies and time delay strategies, can be seen as special forms of non-anticipative strategies \cite{bardi1997optimal}\cite{Cardaliaguet2018}.
	
	Besides the two-side maximum principle and Hamilton-Jacobi-Isaacs equation, the geometric method is another widely used method to solve pursuit-evasion games. The main idea of this method is to take geometric concepts, such as the Apollonius circle, as the research objects to solve the qualitative or quantitative problems in the pursuit-evasion game. Isaacs mentioned this method in his pioneering work, \cite{MR0210469}. He took the Apollonius circle as the safe region of the evader and used this concept to discuss the problem of the interceptor and bomber. Dutkevich and Petrosyan \cite{dutkevich1972games}  perfected the results of Isaacs. Huang \cite{huang2011guaranteed} and Zhou \cite{zhou2016cooperative}  provided a pursuit strategy by Voronoi diagrams in the game with multiple pursuers and a single evader. This strategy can be transformed into a decentralized, real-time algorithm. Chen and Zha used geometric methods to study the qualitative problem of fishing game problems \cite{zha2016construction}, and applied their results in the multi-player pursuit-evasion game with one superior evader \cite{chen2016multi}. Dorothy \cite{dorothy2024one} presented a state-feedback strategy for the pursuer, by which the pursuer can capture the evader in an arbitrarily close neighborhood of the initial Apollonius circle.
	
	Geometric methods are also used in the variations of classic pursuit-evasion games, such as target defense games. Garcia \cite{garcia2018design} \cite{garcia2018optimal} \cite{garcia2019pursuit} and Liang \cite{liang2019differential}\cite{liang2022reconnaissance}\cite{liang2023targets} applied geometric methods in the active target defense game to solve its qualitative and quantitative problems. Yan \cite{yan2019task} proposed a task assignment algorithm for multi-player static target defense games based on the geometric solution of pursuit-evasion games in a convex domain. He also gave the defender's strategy in the 3D multi-player static target defense game by analyzing the generalized concept of Cartesian oval \cite{yan2022matching}. Fu \cite{fu2023justification} and Lee \cite{lee2024solutions} proved that the state feedback strategy provided by the geometric solution satisfies the Hamilton-Jacobi-Isaacs partial differential equation in the multi-defender single-attacker static target defense game, which illustrates the optimality of this strategy from another perspective. Yan \cite{yan2023homicidal}\cite{yan2024multiplayer} designed pursuit strategies of the Dubins car defender by geometric methods. 
	The geometric concepts involved in the above articles are all related to the dominance region of the evader.
	
	In recent years, there have been many studies on pursuit-evasion games in the presence of obstacles. Margellos \cite{margellos2011hamilton} and Fisac \cite{fisac2015reach} presented a method to compute the reachable sets in reach-avoid differential games by the Hamilton-Jacobi-Isaacs equations, which can be used to deal with pursuit-evasion games in the presence of obstacles. However, this method cannot be directly applied in large-scale multi-player pursuit-evasion games because finding the numerical solution of the Hamilton-Jacobi-Isaacs equation suffers from the curse of dimensionality. Bhattachary \cite{bhattacharya2016visibility} and Zou \cite{zou2018optimal} considered a visibility-based target tracking game in an environment containing an obstacle. In this problem, an observer attempts to maintain a persistent line-of-sight with a mobile target. \cite{bhattacharya2016visibility} considered the optimal control policies and the trajectories for the players near the usable part on the terminal manifold in a game with a circular obstacle. The author of \cite{zou2018optimal} summarized the optimal trajectories of the observer at all initial positions in a game with a single corner obstacle. Bhadauria \cite{bhadauria2012capturing} studied a pursuit-evasion game with multiple pursuers and one evader in bounded nonconvex regions. The authors presented a pursuit strategy algorithm that can guarantee the capture of the pursuers and estimated the time complexity of the algorithm.
	There has been limited research on dominance regions in scenarios involving obstacles. Oyler \cite{oyler2016pursuit} studied the dominance regions of pursuit-evasion games in the presence of obstacles. Two methods for constructing the boundary of the dominance region were proposed, and the effects of obstacles that affect the players asymmetrically were discussed. Oyler \cite{oyler2016pursuit} also applied dominance regions in the predator-prey-protector game. However, the author of \cite{oyler2016pursuit} did not emphasize the information pattern of the problems he studied. Based on the overall content of the entire article, we believe that they were considered under the open-loop information pattern. Under this information pattern, the pursuer may know the future behavior of the evader and formulate its current control input by this. It is anticipative. In engineering applications, it is more practical to study the non-anticipative strategy.
	
	This paper mainly studies the properties of the evader's dominance region, especially in the non-anticipative information pattern. 
	The author of \cite{oyler2016pursuit} showed that the initial dominance region of the evader is its reachable region in the Stackelberg open-loop information pattern (though there was not a formal and rigorous proof). The implications of this are two-fold:
	First, the evader can reach any position of its initial dominance region regardless of information patterns because any pursuer's strategy corresponds to an open-loop representation when the evader's control input function is determined \cite{MR1311921}. The pursuer can not prevent the evader leaving any proper subset of its initial dominance region.
	Second, the pursuer can prevent the evader leaving its initial dominance region  if the pursuer knows the evader's control inputs in the whole time domain at the initial moment, which means that the pursuer uses anticipative strategy.  
	Thus, the core problem in this article is under what condition the pursuer has a non-anticipative strategy to prevent the evader leaving its initial dominance region before capture regardless of evader's control input. We call such a strategy PELIDR non-anticipative strategy (where ``PELIDR" means ``Preventing the Evader Leaving its Initial Dominance Region"). 
	The main contributions of this paper are as follows. \\
	(a) We use the shortest path distance to define the dominance region of the evader and rigorously prove that it is the reachable region of the evader in the open-loop sense (\thref{reachable_is_domainance1} and \thref{reachable_is_domainance2}). \\
	(b) We prove that the pursuer's PELIDR non-anticipative strategy exists in the absence of obstacles  (\thref{therom1}). \\
	(c) A counter example is presented to illustrate that the pursuer's PELIDR non-anticipative strategies may not exist in the presence of obstacles (\thref{con_example}), and we provide a necessary condition for the existence of PELIDR non-anticipative strategies in general situations with obstacles (\thref{necessary_condition}). \\
	(d) We provide a sufficient condition for the pursuer to have a PELIDR non-anticipative strategy in the situation where there is only one corner obstacle (\thref{simple_key_theorem}). \\
	(e) We explore the application of the evader's dominance regions in static target defense games under the non-anticipative information pattern.
	
	The remainder of the paper is organized as follows. In Section 2, we introduce some notions and terminologies and show the reachability of the evader's dominance region in the open-loop sense. In Section 3, we consider the dominance regions in the absence of obstacles. In Section 4, we consider the dominance regions in general cases with obstacles. In Section 5, we consider the dominance regions in the case with a single corner obstacle. In Section 6, we apply dominance regions in static target defense games. Finally, we conclude the paper in Section 7.

	\section{Preliminaries}
	\subsection{Notions}
	Given topological space $(X,\Omega)$ and $A \subseteq X$, $A^{\circ}$ denotes the interior of $A$, $\overline{A}$ denotes the closure of $A$, and $\partial A$ denotes the boundary of $A$.\\
	Given $x\in\mathbb{R}^2$, $\left\| x\right\| $ denotes the Euclidean norm of $x$.\\
	Given $x_0\in \mathbb{R}^i,r>0$, $B_i(x_0,r)\triangleq \left\lbrace x\in\mathbb{R}^i:\left\| x-x_0\right\| <r \right\rbrace$.
	\subsection{Cartesian ovals and Apollonius circles}
\begin{figure}[h]
		\centering
		\includegraphics[height=3cm]{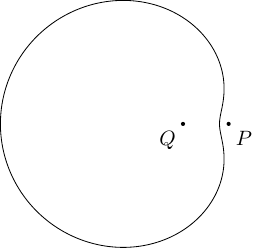}
		\caption{Cartesian oval of the second type.}
		\label{zintro_oval002} 
	\end{figure}

	Let $P$ and $Q$ be fixed points in the plane. $\left| PS\right| $ and $\left| QS\right| $ denote the Euclidean distances from these points to a third variable point $S$ and $m,a\in\mathbb{R}$. Then, the Cartesian oval is the locus of points $S$ satisfying 
	\begin{align*}
		\left| PS\right|+m\left| QS\right| =a,
	\end{align*}
	where $P,Q$ are called the focal points of the Cartesian oval. In this article, there are two special types of Cartesian ovals involved.\\
	\begin{figure}[h]
		\centering
		\includegraphics[height=3cm]{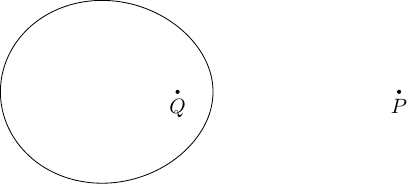}
		\caption{Cartesian oval of the first type.}
		\label{zintro_oval001}
	\end{figure}
	\begin{figure}[h]
		\centering
		\includegraphics[height=3cm]{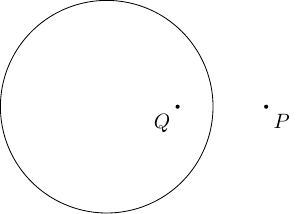} 
		\caption{Apollonius circle.}
		\label{zintro_circle001}
	\end{figure}\\
	For the first type, $m<-1, a\in[0,\left| PQ\right|)$. As Fig.~\ref{zintro_oval001} shows, the Cartesian oval of this type is a bounded,  strictly convex, simple, closed curve \cite{dutkevich1972games}\cite{fu2023justification}\cite{yan2021optimal}. The interior of the Cartesian oval is a bounded strictly convex set. $P$ is outside and $Q$ is inside the Cartesian oval. In this article, we call this type of Cartesian oval the first type of Cartesian oval, with $P$ and $Q $ being the outer and inner focal points respectively. In particular, when $a=0$, the Cartesian oval degenerates into an Apollonius circle, as shown in Fig.~\ref{zintro_circle001}.\\
	For the second type, $m<-1,a\in(m\left| PQ\right|,0)$. As Fig.~\ref{zintro_oval002} shows, the Cartesian oval in this type is bounded, simple, closed but not convex curve. $P$ is outside the Cartesian oval, and $Q$ is inside the Cartesian oval. In this article, we call this type of Cartesian oval the second type of Cartesian oval, with $P$ and $Q$ being the outer and inner focal points, respectively. 
	\subsection{Pursuit-evasion games with or without obstacles}\label{PEgame}
	Consider a pursuit-evasion game taking place in a connected, closed region $X$ in  $\mathbb{R}^2$. There are two players in the game: one pursuer and one evader. Both of them are assumed to be mass points with Isaacs' simple motion\cite{MR0210469}. We also assume that the speeds of the players are constants and the pursuer is faster than the evader. Without loss of generality, we assume that the speed of the evader is $1$ and the speed of the pursuer is $\alpha$ ($\alpha>1$). The players are not allowed to enter $\mathbb{R}^2\setminus X$, which can be considered as obstacles. Let $l\ge0$ be the capture radius. Assume that the distance between the pursuer and the evader at the initial moment is greater than $l$. When the distance between the pursuer and the evader is less than or equal to $l$, the pursuer captures the evader and the game ends. For games with obstacles, we only consider the case where $l=0$. Let $x_p,x_e\in\mathbb{R}^2$ be the location coordinates of the pursuer and the evader, respectively. The equations of motion are as follows:
	\begin{equation}\label{dynamic system} 
		\begin{split}
			&\dot{x}_p(t)=\alpha u_p(t) ,\dot{x}_e (t)= u_e (t),  \\
			&x_p(t_0)=x_{p0},x_e(t_0)=x_{e0}, 
		\end{split}
	\end{equation}
	where $x_{p0},x_{e0}$ represent the initial positions and $u_p(t),u_e(t)$ are the control input functions of the pursuer and the evader, respectively. Assume that the motion trajectories of the pursuer and the evader are continuous and piecewise smooth and the players is allowed to change their movement direction suddenly. Let 
	\begin{align*}
		&\mathcal{U}_{t_0} \triangleq \\
		&\left\lbrace \right. u:\mathbb[t_0,+\infty)\rightarrow \partial B_2(\boldsymbol{0},1):
		\left.u \text{ is piecewise continuous} \right\rbrace. 
	\end{align*}
	Let $\mathcal{U}_{t_0}^p(x_{p0})$ and $\mathcal{U}_{t_0}^e(x_{e0})$ denote the sets of
	 admissible control input functions for the pursuer and the evader respectively. Because the players are not allowed to enter $\mathbb{R}^2\setminus X$, $\mathcal{U}_{t_0}^p(x_{p0})$ or $\mathcal{U}_{t_0}^e(x_{e0})$ may not be $\mathcal{U}_{t_0}$, but they are subsets of $\mathcal{U}_{t_0}$.
	For any control input functions $u_p(\cdot),u_e(\cdot)$ , any initial time $t_0$, and any initial state $\mathbf{x}_0=(x_{p0},x_{e0})$, the solution of \eqref{dynamic system} is denoted by $\mathbf{x}(t;t_0,\mathbf{x}_0,u_p,u_e)=(x_p(t;t_0,x_{p0},u_p),x_e(t;t_0,x_{e0},u_e))$. For simplicity, it can be abbreviated as $\mathbf{x}(t;u_p,u_e)=(x_p(t;u_p),x_e(t;u_e))$ or $\mathbf{x}(t)=(x_p(t),x_e(t))$ when there is no ambiguity. The terminal set of the game is
	\begin{align*}
		\mathcal{T}=\left\lbrace (x_p,x_e):\left\| x_p-x_e\right\| \le l \right\rbrace.  
	\end{align*}
	The terminal time of the game is
	\begin{align}
		t_f(t_0,\mathbf{x}_0,u_p,u_e)=\inf \left\lbrace t :\mathbf{x}(t;t_0,\mathbf{x}_0,u_p,u_e) \in \mathcal{T}  \right\rbrace \label{terminal_time}.
	\end{align}
	This can be abbreviated as $t_f(u_p,u_e)$ or $t_f$ when there is no ambiguity. We define the infimum of the empty set as $+\infty$. $t_f=+\infty$ means that the pursuer can not capture the evader with the control input functions.
	\subsection{Information pattern}
	Information patterns \cite{MR1311921}\cite{olsder1977information} play an important role in differential games. In \cite{oyler2016pursuit}, the author did not emphasize the information pattern of the problems he studied. Based on the overall content of the entire article, we believe that they were considered under the Stackelberg open-loop information pattern. Under this information pattern, the pursuer may know the future behavior of the evader and formulate its current control input by this. In this article, we adopt the non-anticipative information pattern. The concept of ``non-anticipative strategy'' was proposed over 50 years ago \cite{varaiya1967existence}\cite{roxin1969axiomatic}\cite{elliott1972existence}, and it was widely used in several studies \cite{evans1984differential}\cite{soravia1993pursuit}\cite{cardaliaguet1996differential}\cite{mitchell2005time}\cite{margellos2011hamilton}\cite{fisac2015reach}. Non-anticipative strategies do not depend on the future behavior of the other player. The executor of a non-anticipative strategy can formulate its current control input only by the state and the control input of the other players from the initial moment to the current moment. For example, state feedback strategies and time-delay strategies both can be seen as non-anticipative strategies \cite{bardi1997optimal}\cite{Cardaliaguet2018}. 
	
	In this paper, a strategy for the pursuer is a map, from $\mathcal{U}_{t_0}^e(x_{e0})$ to $\mathcal{U}_{t_0}^p(x_{p0})$ . This means that the pursuer answers to each control input function $u_e\in\mathcal{U}_{t_0}^e(x_{e0})$ of the evader by a control input function $u_p=\delta[u_e]\in \mathcal{U}_{t_0}^p(x_{p0})$. However, since we wish to formalize the fact that the pursuer can not guess in advance the future behavior of the evader, we require that such a map is non-anticipative.
	\begin{defn}
		Given initial moment $t_0$ and initial state $\mathbf{x}_0$, a map $\delta$ from $\mathcal{U}_{t_0}^e(x_{e0})$ to $\mathcal{U}_{t_0}^p(x_{p0})$ is called the pursuer's non-anticipative strategy, if it is non-anticipative:
		\begin{equation}\label{nonanticipate}
			\begin{split}
				&u_e(t)=u_e'(t), \text{ a.e.}[t_0,t'] \\
				\Rightarrow &\delta[u_e](t)=\delta[u_e'](t), \text{ a.e.}[t_0,t'].
			\end{split}
		\end{equation}
		The set of all non-anticipative strategies of the pursuer is denoted by 
		\begin{align*}
			\Delta^{t_0}(\mathbf{x}_0)
			=\left\lbrace \delta:\mathcal{U}_{t_0}^e(x_{e0})\rightarrow\mathcal{U}_{t_0}^p(x_{p0}):\right.\\
			\forall t'>t_0,
			u_e(t)=u_e'(t), \text{ a.e.}[t_0,t'] \\
			\left.\Rightarrow \delta[u_e](t)=\delta[u_e'](t), \text{ a.e.}[t_0,t'] \right\rbrace. 
		\end{align*}
	\end{defn}
	Now, we give a way to generate non-anticipative strategies of the pursuer. Given initial moment $t_0$, initial state $\mathbf{x}_0$ and the evader's control input function $u_e(\cdot)$, the control input of the pursuer at time $t$ is only related to the state and the control input of the evader at time $t$:
	\begin{equation}\label{nonanti_stra1}
		u_p(t)=\gamma(x_p(t),x_e(t),u_e(t)). 
	\end{equation}
	If the map $\gamma: \mathbb{R}^2 \times \mathbb{R}^2 \times \partial B_2(\boldsymbol{0},1) \rightarrow \mathbb{R}^2 $ is continuous and locally Lipschitz continuous with respect to $x_p(t)$ for any $x_e(t),u_e(t)$, then there exists a unique solution of
	\begin{align*}
		&\dot{x}_p(t)=\alpha \gamma(x_p(t),x_e(t,u_e),u_e(t)),\\
		&x_p(t_0)=x_{p0},
	\end{align*}
	which is denoted by $\xi_p(t;\gamma,u_e)$. We can generate a map $\delta_p$ from $\mathcal{U}^e_{t_0}(x_{p0})$ to $\mathcal{U}^p_{t_0}(x_{e0})$, which is defined as follows:
	\begin{equation}\label{nonanti_stra2}
		\begin{split}
			\delta_p[u_e](t)=\gamma(\xi_p(t;\gamma,u_e),x_e(t;u_e),u_e(t)). 
		\end{split}
	\end{equation}
	It is obvious that $\delta_p[\cdot]$ is non-anticipative.

	The above method for generating the pursuer's non-anticipative strategies requires the pursuer to know the evader's current control input. In this game, the states are positions of the players, and the control inputs are the velocity directions of the players. It is realistic to assume that the pursuer knows the positions of all players from initial moment to current moment. Subsequently, the pursuer knows the left derivative of evader's position with respect to time at current time:
	\begin{align*}
		\lim_{\Delta t\rightarrow 0^+}\frac{x_e(t)-x_e(t-\Delta t)}{\Delta t}.
	\end{align*}
	We also assume that the control input function is piecewise continuous. Subsequently, the left derivative at current time is equal to the current control input almost everywhere:
	\begin{align*}
		\lim_{\Delta t\rightarrow 0^+}\frac{x_e(t)-x_e(t-\Delta t)}{\Delta t}=u_e(t), \:\text{a.e. } t\in(t_0,+\infty).
	\end{align*}
	Thus, the pursuer knows the evader's current control input in the time domain almost everywhere. The integrals of two functions that are equal almost everywhere are equal. Therefore, we believe it is realistic that the pursuer has access to the
	evader's current velocity $u_e(t)$, theoretically.
	In terms of engineering application, there are several methods of estimating the velocity of opponent moving object \cite{battistini2014differential}\cite{fonod2017blinding}\cite{mohanan2020target}. In addition, the velocity of a moving object can be accurately measured by advanced sensor systems  , such as phased-array radars \cite{pihl2012phased}.\\
	For practical applications, we can also approximate the evader's current velocity by the delayed availability, i.e. replacing $u_e(t)$ by $u_e(t-\Delta t)$ (where $\Delta t >0$) in \eqref{nonanti_stra1}. The evaluation of this approximation deserves further study. We leave it in future works.\\
	If the state is on the dispersal surface \cite{MR0210469}, there are some challenges without the access to the opponent's current control input. As shown in Appendix \ref{challenge}, one player is free to choose between two equally optimal alternatives while its opponent can only behave optimally if it knows which choice the player makes. Some researchers have addressed the dilemma caused by the singularity in some special cases, such as \cite{milutinovic2021rate} and \cite{dorothy2024one}.\\
	Even if the assumption of knowing the opponent's current control input is unrealistic, our study is still meaningful. In \cite[Chapter VIII]{bardi1997optimal}, the authors demonstrated that the value function of non-anticipative information pattern with the access to the opponent's current control input is equal to the value function of feedback information pattern without the access. This implies that strategies with the access and strategies without the access have the same performance limit. Studying non-anticipative strategies with the access helps us to explore the performance limit of feedback strategies without the access.

	\subsection{Shortest path distance and dominance region}
	We define a binary function on $X$. For any $x_1,x_2\in X$,
	\begin{align*}
		d_L(x_1,x_2)\triangleq \inf \left\lbrace L(\mathbf{r}):\mathbf{r}:[a,b]\rightarrow X ,\mathbf{r}(a)=x_1,\right.\\
		\left. \mathbf{r}(b)=x_2,\mathbf{r}\text{ is continuous and piecewise smooth}\right\rbrace,
	\end{align*}
	where $L(\mathbf{r})$ represents the length of the curve $\mathbf{r}$. $d_L$ is a metric on $X$ \cite[Example 2.2.2]{MR1835418}. We call this the shortest path distance on $X$ in this article. $d_L$ is Euclidean distance when $X$ is convex. $X$ is closed, so for any $x_1,x_2\in X$, there exist a piecewise smooth curve $\mathbf{r}$ connecting $x_1,x_2$ such that $L(\mathbf{r})=d_L(x_1,x_2)$ \cite[Theorem 2.5.23]{MR1835418}. We call this curve $\mathbf{r}$ the shortest path connecting $x_1$ and $x_2$. If the shortest path connecting $x_1$ and $x_2$ is the line segment connecting $x_1$ and $x_2$, we say that $x_1$ and $x_2$ are visible to each other. For $x\in X$, the set of all visible points in $X$ to $x$ is called the visible region of $x$, denoted by $\mathcal{V}(x)$. 
	\begin{defn}
		In the case without obstacles, we allow the capture radius to be a positive number or zero. For any $x_1,x_2\in\mathbb{R}^2$, let
		\begin{align*}
			\mathcal{D}_l(x_1,x_2)\triangleq \left\lbrace x\in\mathbb{R}^2 : \left\| x-x_1\right\|-\alpha \left\| x-x_2 \right\| > l  \right\rbrace.
		\end{align*}
		When $\left\| x_1-x_2\right\|>l $, the boundary of $\mathcal{D}_l(x_1,x_2)$ is the Cartesian oval of the first type with $x_1,x_2$ being the focal points. Given the positions of the pursuer and evader $x_p,x_e$, in the case without obstacles, we call $\mathcal{D}_l(x_p,x_e)$  the evader's dominance region.\\
		For the games with obstacles, we only consider the case where $l=0$. For any $x_1,x_2\in X$, let
		\begin{align*}
			\mathcal{D}(x_1,x_2) \triangleq \left\lbrace x\in X : d_L(x,x_1)-\alpha d_L(x,x_2) > 0  \right\rbrace.
		\end{align*}
		Given the positions of the pursuer and evader $x_p,x_e$, in the case with obstacles, we call $\mathcal{D}(x_p,x_e)$ the evader's dominance region.
	\end{defn}
	\subsection{Reachable region of the evader}
	\begin{defn}
		Given initial time $t_0$ and initial state $\mathbf{x}_0=(x_{p0},x_{e0})$, the reachable region of the evader is defined as follows:
		\begin{align*}
			\mathcal{R}_e(t_0,\mathbf{x}_0) \triangleq \left\lbrace x\in X :\exists u_e\in \mathcal{U}_{t_0}^e(x_{e0}),\forall u_p\in\mathcal{U}_{t_0}^p(x_{p0}),\right.\\
			\left.\exists t\in [t_0,t_f(u_p,u_e)) ,x_e(t;t_0,x_{e0},u_e)=x\right\rbrace. 
		\end{align*}
		The set defined above is independent of $t_0$ due to the time-invariance of \eqref{dynamic system}, thus $\mathcal{R}_e(t_0,\mathbf{x}_0) $ can be abbreviated as $\mathcal{R}_e(\mathbf{x}_0)$ or $\mathcal{R}_e(x_{p0},x_{e0})$. 
	\end{defn}
	Next, we show the relationship between the dominance region and the reachable region of the evader.
	\begin{thm}\thlabel{reachable_is_domainance1}
		In the case without obstacles, 
		\begin{align*}
			\mathcal{R}_e(x_{p0},x_{e0})=\mathcal{D}_l(x_{p0},x_{e0}).
		\end{align*}
	\end{thm}
	\begin{thm}\thlabel{reachable_is_domainance2}
		In the case with zero capture radius and obstacles, 
		\begin{align*}
			\mathcal{R}_e(x_{p0},x_{e0})=\mathcal{D}(x_{p0},x_{e0}). 
		\end{align*}
	\end{thm}
	It is easy to prove \thref{reachable_is_domainance1} by the convexity of Cartesian ovals of the first type. A brief proof of \thref{reachable_is_domainance1} can be found in \cite{yan2021optimal}. Next, we present a proof of \thref{reachable_is_domainance2}. We first prove a lemma. 
	\begin{lem}\thlabel{shortestline_in_domainance}
		For any $x_{p0},x_{e0}\in X$, $(x_{p0}\neq x_{e0})$ and $x \in \overline{\mathcal{D}(x_{p0},x_{e0})}$, let $\mathbf{r}^*_x(\cdot)$ be a shortest path connecting $x_{e0}$ and $x$ with $\mathbf{r}^*_x(a)=x_{e0},\mathbf{r}^*_x(b)=x$. Then $\left\lbrace \mathbf{r}^*_x(s) : s \in [a,b) \right\rbrace \subseteq \mathcal{D}(x_{p0},x_{e0}) $. 
	\end{lem}
	\begin{pf}
		For the sake of contradiction, assume that $\exists s \in [a,b),r^*_x(s) \notin \mathcal{D} (x_{p0},x_{e0})$, i.e. $d_L( r^*_x(s),x_{p0}) \le \alpha d_L ( r^*_x(s),x_{e0})$. By the triangle inequality of metric $d_L$ \cite{MR1835418} and $\alpha>1$, 
		\begin{align*}
			d_L(x,x_{p0}) &\le d_L( r^*_x(s),x_{p0}) + d_L( r^*_x(s),x) \\
			&< d_L( r^*_x(s),x_{p0}) + \alpha d_L( r^*_x(s),x) \\ 
			&\le \alpha d_L( r^*_x(s),x_{e0}) + \alpha d_L( r^*_x(s),x) \\
			&= \alpha d_L(x,x_{e0}).       		
		\end{align*}
		The above inequality contradicts that $x \in \overline{\mathcal{D}(x_{p0},x_{e0})}$. The proof is completed. \qed
	\end{pf}
	\begin{pf}[Proof of \thref{reachable_is_domainance2}]
		Let
		\begin{align*}
			&t_p(x;t_0,x_{p0})\\
			=&\inf\left\lbrace t:\exists u_p\in \mathcal{U}_{t_0}^p(x_{p0}), x_p(t;t_0,x_{p0},u_p)=x \right\rbrace, \\
			&t_e(x;t_0,x_{e0})\\
			=&\inf\left\lbrace t:\exists u_e\in \mathcal{U}_{t_0}^e(x_{e0}), x_e(t;t_0,x_{e0},u_e)=x \right\rbrace.
		\end{align*}
		$t_p(x;t_0,x_{p0})$ is the minimum moment when the pursuer reaches $x$ under the initial condition $t_0,x_{p0}$. $t_e(x;t_0,x_{e0})$ is the minimum moment when the evader reaches $x$ under the initial condition $t_0,x_{e0}$. It is obvious that
		\begin{align*}
			&t_p(x;t_0,x_{p0})=t_0+\frac{1}{\alpha}d_L(x_{p0},x),\\
			&t_e(x;t_0,x_{e0})=t_0+d_L(x_{e0},x).
		\end{align*}
		$\mathcal{R}_e(x_{p0},x_{e0})=\mathcal{D}(x_{p0},x_{e0})$ is equivalent to $\mathcal{D}(x_{p0},x_{e0}) \subseteq \mathcal{R}_e(x_{p0},x_{e0})$ and $\mathcal{D}(x_{p0},x_{e0})^c \subseteq \mathcal{R}_e(x_{p0},x_{e0})^c$.
		\par
		We first prove that $\mathcal{D}(x_{p0},x_{e0}) \subseteq \mathcal{R}_e(x_{p0},x_{e0})$. Pick $x \in \mathcal{D}(x_{p0},x_{e0})$. Let $r^*_x(\cdot)$ be a shortest path connecting $x_{e0}$ and $x$, with $r^*_x(a)=x_{e0},r^*_x(b)=x$. Let the evader select a control input function $u_e(\cdot)\in \mathcal{U}_e^{t_0}(x_{e0})$ such that the evader goes to $x$ along $r^*_x(\cdot)$ under the initial conditions $t_0$ and $x_{e0}$. By \thref{shortestline_in_domainance}, the shortest path is in $\mathcal{D}(x_{p0},x_{e0})$. Then
		\begin{align*}
			t_p(r^*_x(s),t_0,x_{p0})>t_e(r^*_x(s),t_0,x_{e0}),\forall s \in [a,b].
		\end{align*}
		Thus, the pursuer cannot reach any point of $r^*_x(\cdot)$ at the same time when the evader reaches it by $u_e(\cdot)$. The evader can reach $x$ before being captured. We get $\mathcal{D}(x_{p0},x_{e0}) \subseteq \mathcal{R}_e(x_{p0},x_{e0})$. 
		\par
		Next, we prove that $\mathcal{D}(x_{p0},x_{e0})^c \subseteq \mathcal{R}_e(x_{p0},x_{e0})^c$. Pick $x \in \mathcal{D}(x_{p0},x_{e0})^c$. Let the evader select a control input function $\bar{u}_e(\cdot)\in \mathcal{U}_e^{t_0}(x_{e0})$ such that the evader passes $x$, i.e. $\exists t \in (t_0,+\infty),x(t;t_0,x_{e0},\bar{u}_e)=x$. Let
		\begin{align*}
			\tau^*=&\tau^*(t_0,x_{p0},x_{e0},\bar{u}_e)\\
			\triangleq&
			\inf\left\lbrace t:x_e(t;t_0,x_{e0},\bar{u}_e)\in\mathcal{D}(x_{p0},x_{e0})^c\right\rbrace. 
		\end{align*}
		$\tau^*$ is the moment when the evader first leaves $\mathcal{D}(x_{p0},x_{e0})$ under the control input function $\bar{u}_e(\cdot)$. Let $x^*\triangleq x_e(\tau^*;t_0,x_{e0},\bar{u}_e)$. Then, we have
		\begin{align*}
			\tau^*&\ge \inf\left\lbrace t:\exists u_e\in \mathcal{U}_{t_0}^e(x_{e0}), x_e(t;t_0,x_{e0},u_e)=x^* \right\rbrace\\
			&=t_e(x^*;t_0,x_{e0})
		\end{align*}
		It is obvious that $x^*\in\partial \mathcal{D}(x_{p0},x_{e0})$. Let the pursuer go to $x^*$ along the shortest path connecting $x_{p0}$ and $x_e(\tau^*;t_0,x_{e0},\bar{u}_e)$. The moment when the pursuer reaches $x^*$ is $t_p(x^*;t_0,x_{p0})$. By the definition of $\mathcal{D}(x_{p0},x_{e0})$,
		\begin{align*}
			t_p(x^*;t_0,x_{p0})=t_e(x^*;t_0,x_{e0}).
		\end{align*}
		Without loss of generality, assume that the pursuer does not meet the evader in the time interval $[t_0,t_p(x^*;t_0,x_{p0}))$. 
		If $t_e(x^*;t_0,x_{e0})= \tau^*$, then $t_f=\tau^*$. The evader is always in $\mathcal{D}(x_{p0},x_{e0})$ during $[t_0,t_f)$ so it does not reach $x$ before being captured.
		If $t_e(x^*;t_0,x_{e0})< \tau^*$, let the pursuer retrace the trajectory of the evader from the moment $t_e(x^*;t_0,x_{e0})$ and the position $x^*$ untill the pursuer meets the evader and captures it. The moment when they meets is simply $t_f$. It is obvious that $t_f<\tau^*$. Thus, the evader is always in $\mathcal{D}(x_{p0},x_{e0})$ during $[t_0,t_f)$, and the evader cannot reach $x$ before being captured. We get $\mathcal{D}(x_{p0},x_{e0})^c \subseteq \mathcal{R}_e(x_{p0},x_{e0})^c$ .\qed
	\end{pf}
	The complement of $\mathcal{R}_e(\mathbf{x}_0)$ is represented as follows:
	\begin{equation}\label{complement of reachable set}
		\begin{split}
			&\mathcal{R}_e(\mathbf{x}_0)^c \\
			= &\left\lbrace x\in \mathbb{R}^2:\right.\left.\forall u_e\in \mathcal{U}^e_{t_0}(x_{e0}), \exists u_p\in\mathcal{U}^p_{t_0}(x_{p0}),\right.\\
			&\left.\forall t\in [t_0,t_f(u_p,u_e)) ,x_e(t;u_e,x_{e0}) \neq x\right\rbrace. 
		\end{split}
	\end{equation}
	In other words, for any point in the complement of the evader's initial dominance region, for any control input function of the evader, there exists a control input function of the pursuer by which the pursuer can prevent the evader from reaching the point before being captured. It should be noted that the control input function of the pursuer in the previous sentence may depend on all the evader's control input in the time interval $[t_0,+\infty)$ according to \eqref{complement of reachable set}; this is an anticipative strategy. The game can be seen as an open-loop Stackelberg game \cite{MR1311921}. We naturally raise the question of whether the pursuer can prevent the evader from reaching the complement of the evader's initial dominance region before being captured in the non-anticipative information pattern i.e. whether the pursuer has a PELIDR non-anticipative strategy. We investigate this issue in the next three sections.

	\section{Existence of a PELIDR non-anticipative strategy without obstacles}\label{section without obstacle}
\begin{figure}[h]
		\centering
		\includegraphics[height=3cm]{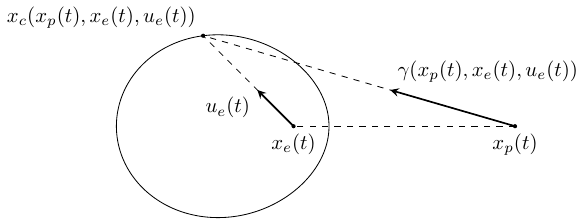}
		\caption{The non-anticipative strategy $\delta_{t_0,\mathbf{x}_0}^*$ is generated by $\gamma$ like \eqref{nonanti_stra2} with initial time $t_0$ and initial state $\mathbf{x}_0=(x_{p0},x_{e0})$.}
		\label{znonanticipate_strategy} 
	\end{figure}

	In this section, we prove that there exists a PELIDR non-anticipative strategy of the pursuer in the absence of obstacles (i.e. $X=\mathbb{R}^2$).\\
	We first define the PELIDR non-anticipative strategy in the form of \eqref{nonanti_stra1}. 
	Given $x_p(t),x_e(t),u_e(t)$, draw a ray from $x_e(t)$ along $u_e(t)$. $\partial\mathcal{D}_l(x_p(t),x_e(t))$ is a strictly convex, simple, closed curve, and $x_e(t)$ is in the interior of $\partial\mathcal{D}_l(x_p(t),x_e(t))$. Thus, there exists a unique intersection point of the ray and $\partial\mathcal{D}_l(x_p(t),x_e(t))$. The intersection point is related to $x_p(t),x_e(t)$ and $u_e(t)$. Let $x_c(x_p(t),x_e(t),u_e(t))$ denote the coordinate of the intersection point. The control input of the pursuer at time $t$ is defined as follows:
	\begin{equation*}
		\begin{split}
			u_p(t)&=\gamma(x_p(t),x_e(t),u_e(t))\\
			&=\frac{x_c(x_p(t),x_e(t),u_e(t))-x_p(t)}{\left\| x_c(x_p(t),x_e(t),u_e(t))-x_p(t)\right\| },
		\end{split}
	\end{equation*}
	as shown in Fig.~\ref{znonanticipate_strategy}. 
	By extending this function properly, it can generate a non-anticipative strategy like \eqref{nonanti_stra2} with initial time $t_0$ and initial state $\mathbf{x}_0=(x_{p0},x_{e0})$. This non-anticipative strategy is denoted as $\delta^*_{t_0,\mathbf{x}_0}$, abbreviated as $\delta^*$.
	\begin{thm}\thlabel{therom1}
		In the absence of obstacles, given initial time $t_0$ and initial state $\mathbf{x}_0=(x_{p0},x_{e0})$,\\
		(a)there exists a constant $c>0$ such that $t_f(\delta^*[u_e],u_e)<c$ for any $u_e(\cdot)\in \mathcal{U}_{t_0}$;\\
		(b)for any $u_e\in \mathcal{U}_{t_0}$, for any $t\in [t_0,t_f(\delta^*[u_e],u_e))$, $x_e(t;u_e)\in \mathcal{D}_l(x_{p0},x_{e0})$.
	\end{thm}
	The authors of \cite{dutkevich1972games} proved that if the evader moves along a straight line, the pursuer can capture the evader in $\overline{\mathcal{D}_l(x_{p0},x_{e0})}$ by the non-anticipative strategy mentioned above. We provide a more general proof by a lemma as follows.
	
	\begin{lem}\thlabel{lemma1}
		$x_p,x_e\in \mathbb{R}^2,\left\| x_p-x_e\right\|  > l$. For any $x_1,x_2\in \partial \mathcal{D}_l(x_p,x_e)$,
		\begin{align*}
			\frac{(x_1-x_p)\cdot(x_2-x_p)}{\left\| x_1-x_p\right\| \left\| x_2-x_p\right\| }\ge \frac{(x_1-x_e)\cdot(x_2-x_e)}{\left\| x_1-x_e\right\| \left\| x_2-x_e\right\| }.
		\end{align*}
	\end{lem}
	\begin{figure}[h]
		\centering
		\subfigure[$\chi(x)<\pi$]{
			\includegraphics[height=3cm]{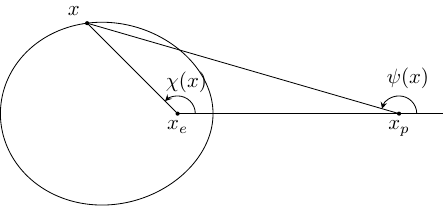}
		}
		\subfigure[$\chi(x)>\pi$]{
			\includegraphics[height=3cm]{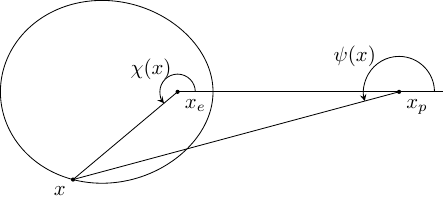}
		}
		\caption{$\chi(x)$ and $\psi(x)$ in the proof of \thref{lemma1}.}
		\label{figure_chi_psi_1}
	\end{figure}
	\begin{pf}
		Referring to Fig.~\ref{figure_chi_psi_1}, for any $x\in \partial \mathcal{D}_l(x_p,x_e)$, let $\chi=\chi(x)$ denote the angle between vectors $x-x_e$ and  $x_p-x_e$ ($x_p-x_e$ is the starting edge, $x-x_e$ is the ending edge, and counterclockwise is the positive direction). Let $\psi=\psi(x)$ denote the angle between vectors $x-x_p$ and  $x_p-x_e$ ($x_p-x_e$ is the starting edge, $x-x_p$ is the ending edge, and counterclockwise is the positive direction). Without loss of generality, assume that $\chi(x_2)-\chi(x_1)\in[0,\pi]$. If
		\begin{align}\label{c-oval_local_pro}
			\left| \frac{ d\psi(x)}{d\chi(x)}\right| \le 1
		\end{align}
		holds on the whole closed curve $\partial \mathcal{D}_l(x_p,x_e)$, we obtain
		\begin{align*}
			&\left| \psi(x_2)-\psi(x_1)\right| =\left| \int_{\chi(x_1)}^{\chi(x_2)} \frac{ d\psi(x)}{d\chi(x)} d\chi(x)\right| \\
			\le& \int_{\chi(x_1)}^{\chi(x_2)} \left| \frac{ d\psi(x)}{d\chi(x)} \right|d\chi(x)  \le \int_{\chi(x_1)}^{\chi(x_2)} d\chi(x) \\
			\le& \chi(x_2)-\chi(x_1).
		\end{align*}
		Then,
		\begin{align*}
			\cos(\psi(x_2)-\psi(x_1)) \ge \cos(\chi(x_2)-\chi(x_1)).
		\end{align*}
		The lemma holds by the geometric interpretation of the vectors' inner product. Let us prove that \eqref{c-oval_local_pro} holds on the whole of $\partial \mathcal{D}_l(x_p,x_e)$. Let
		$r_1(x)=\left\| x-x_p\right\|,
		r_2(x)=\left\| x-x_e\right\|,
		\sigma=\left\| x_p-x_e\right\|$.
		$r_1=\alpha r_2+l$ holds. It is only necessary to consider the case $\chi\in[0,\pi]$ due to the symmetry of $\mathcal{D}_l$. By the parametric representation of the Cartesian oval \cite[Lemma 1]{garcia2021cooperative2}, we obtain that $\chi$ and $r_2$ (or $r_1$) correspond one-to-one when $\chi\in[0,\pi]$. The value range of $r_2$ is $[\frac{\sigma-l}{\alpha+1},\frac{\sigma-l}{\alpha-1}]$. By the law of cosines, 
		\begin{align*}
			\cos\chi=\frac{r_2^2+\sigma^2-r_1^2}{2r_2\sigma},\\
			\cos\psi=\frac{r_1^2+\sigma^2-r_2^2}{2r_1\sigma}.
		\end{align*}
		Treating $\cos\chi,\cos\psi, r_1$ as functions of $r_2$, we take the derivative of the above equations on both sides with respect to $r_2$:
		\begin{align*}
			\frac{d\cos\chi}{d r_2}&=\frac{r_1^2+r_2^2-2\alpha r_1 r_2-\sigma^2}{2\sigma r_2^2},\\
			\frac{d\cos\psi}{d r_2}&=\frac{\alpha r_1^2+\alpha r_2^2-2 r_1 r_2-\alpha \sigma^2}{2\sigma r_1^2}.
		\end{align*}
		By the above equations and the law of sines,
		\begin{align*}
			\frac{d\psi}{d\chi}=&\frac{\sin\chi}{\sin\psi} \cdot \frac{d\cos\psi}{d\cos\chi}\\
			=&\frac{r_1}{r_2}\cdot \frac{d\cos \psi /dr_2}{d\cos \chi /dr_2}\\
			=&\frac{\alpha r_1^2r_2-2r_1r_2^2+\alpha r_2^3-\alpha\sigma^2r_2}{r_1^3-2\alpha r_1^2r_2+r_1r_2^2-\sigma^2r_1}.
		\end{align*}
		For convenience, let
		\begin{align}
			f_1&=\alpha r_1^2r_2-2r_1r_2^2+\alpha r_2^3-\alpha\sigma^2r_2,\\
			f_2&=r_1^3-2\alpha r_1^2r_2+r_1r_2^2-\sigma^2r_1 \label{f1}.
		\end{align}
		Substitute $r_2=\frac{1}{\alpha}r_1-\frac{l}{\alpha}$ into \eqref{f1},
		\begin{align*}
			f_2&=r_1\left[ \left(\frac{1}{\alpha^2}-1\right)r_1^2+2l\left(1-\frac{1}{\alpha^2}\right)r_1+\frac{l^2}{\alpha^2}-\sigma^2\right] \\
			&=-\frac{r_1}{\alpha^2}\left[ \left( \alpha^2-1\right) \left( r_1-l\right) ^2 +\alpha^2 \left( \sigma ^2- l^2 \right) \right]. 
		\end{align*}
		Due to $\alpha >1,\sigma>l$, we obtain $\forall r_1>0,f_2(r_1)<0$. Then,
		\begin{align*}
			\left| \frac{f_1}{f_2}\right|\le 1,\forall r_2\in \left[ \frac{\sigma-l}{\alpha+1},\frac{\sigma-l}{\alpha-1}\right]
		\end{align*}
		is converted into
		\begin{align}\label{converted}
			f_1-f_2\ge 0 \bigwedge f_1+f_2\le 0 ,\forall r_2\in \left[ \frac{\sigma-l}{\alpha+1},\frac{\sigma-l}{\alpha-1}\right].
		\end{align}
		Let $g_1=f_1+f_2,\:g_2=f_1-f_2$. Substitute $r_1=\alpha r_2+l$ into $g_1,g_2$ and simplify them:
		\begin{align*}
			g_1(r_2)=l(\alpha^2-1)r_2^2-2\alpha (\sigma^2-l^2)r_2-l(\sigma^2-l^2),\\
			g_2(r_2)=2(\alpha^3-\alpha)r_2^3+3l(\alpha^2-1)r_2^2+l(\sigma^2-l^2).
		\end{align*}
		Consider $g_1$,
		\begin{align*}
			&g_1\left( \frac{\sigma-l}{\alpha+1}\right)\\ =&\frac{-((\alpha+1)l+2\alpha\sigma)(\sigma-l)-l(\alpha+1)(\sigma+l)}{(\alpha+1)/(\sigma-l)}<0,\\
			&g_1\left( \frac{\sigma-l}{\alpha-1}\right)\\ =&\frac{-((\alpha-1)l+2\alpha\sigma)(\sigma-l)-l(\alpha-1)(\sigma+l)}{(\alpha-1)/(\sigma-l)}<0.
		\end{align*}
		$g_1(r_2)$ is a quadratic function with positive quadratic coefficient. Thus, it is a convex function. Then, $g_1(r_2)<0,\forall r_2\in \left[ \frac{\sigma-l}{\alpha+1},\frac{\sigma-l}{\alpha-1}\right]$. The cubic coefficient of $g_2(r_2)$ is positive, and the other coefficients of $g_2(r_2)$ are non-negative. Thus,  $g_2(r_2)>0,\forall r_2\in \left[ \frac{\sigma-l}{\alpha+1},\frac{\sigma-l}{\alpha-1}\right]$. We obtain that \eqref{converted} holds. The proof is completed. \qed
	\end{pf}
	
	Next, we present the proof of \thref{therom1}. 
	\begin{figure}[h]
		\centering
		\includegraphics[height=3cm]{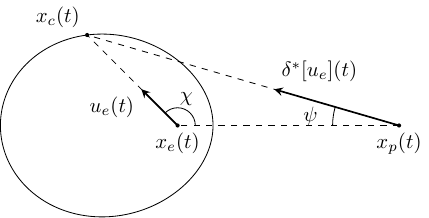} 
		\caption{$\chi(x)$ and $\psi(x)$ in the proof of \thref{therom1}(a).}
		\label{zc_oval003}
	\end{figure}
	\begin{pf}[\thref{therom1}]
		(a) We pick $u_e(\cdot)$ from $\mathcal{U}^e_{t_0}$ and abbreviate $x_c(x_p(t;\delta^*[u_e]),x_e(t;u_p),u_e(t))$ as $x_c(t)$. For $t\in[t_0,t_f(\delta^*[u_e],u_e))$, let $\chi$ denote the angle between vectors $x_c(t)-x_e(t;u_e)$ and $x_p(t;\delta^*[u_e])-x_e(t;u_e)$, and let $\psi$ denote the angle between vectors $x_c(t)-x_p(t;\delta^*[u_e])$ and $x_e(t;u_e)-x_p(t;\delta^*[u_e])$, $\chi ,\psi \in [0,\pi]$, as shown in Fig.~\ref{zc_oval003}. By the law of sines, 
		\begin{align*}
			\frac{\sin\chi}{\sin\psi}
			&=\frac{\left\| x_c(t)-x_p(t;\delta^*[u_e])\right\| }{\left\| x_c(t)-x_e(t;u_e)\right\|}\\
			&=\frac{\alpha\left\| x_c(t)-x_e(t;u_e)\right\| +l}{\left\| x_c(t)-x_e(t;u_e)\right\|}
			\ge\alpha.
		\end{align*}
		Then, $\cos\psi\ge\sqrt{1-\frac{1}{\alpha^2}\sin^2\chi}$. We obtain that
		\begin{align*}
			&\frac{d}{dt}\left\| x_p(t;\delta^*[u_e])-x_e(t;u_e)\right\|  \\
			=&  -\alpha\cos\psi-\cos\chi \\
			\le& -\alpha\sqrt{1-\frac{1}{\alpha^2}\sin^2\chi}-\cos\chi.
		\end{align*}
		Let $f(\chi)=-\alpha\sqrt{1-\frac{1}{\alpha^2}\sin^2\chi}-\cos\chi$,
		\begin{align*}
			f'(\chi)&=\left( \frac{\cos\chi}{\sqrt{\alpha^2-\sin^2\chi}}+1\right) \sin\chi \\
			&>  \frac{\cos\chi+\left| \cos\chi \right|}{\sqrt{\alpha^2-\sin^2\chi} } \sin\chi \ge 0.
		\end{align*}
		Then,
		\begin{align*}
			\frac{d}{dt}\left\| x_p(t;\delta^*[u_e])-x_e(t;u_e)\right\| &\le f(\pi)=1-\alpha,\\
			t_f(\delta^*[u_e],u_e)&\le t_0+\frac{\left\| x_{p0}-x_{e0}\right\|-l}{\alpha-1}.
		\end{align*} 
		(b) Because $x_e(t;u_e)$ is always in $\mathcal{D}_l(x_p(t;\delta^*[u_e]),x_e(t;u_e))$, it is only necessary to prove that $\forall u_e\in \mathcal{U}, \forall t\in [t_0,t_f(\delta^*[u_e],u_e))$, $\overline{ \mathcal{D}_l(x_p(t;\delta^*[u_e]),x_e(t;u_e))} \subseteq \overline{\mathcal{D}_l (x_{p0},x_{e0})}$. Let 
		\begin{equation}\label{level_set}
			\begin{split}
				&\phi(x,t)\\
				=&\left\| x-x_p(t;\delta^*[u_e])\right\|-\alpha \left\| x-x_e(t;u_e) \right\| - l.
			\end{split}
		\end{equation}
		\begin{figure}[h]
			\centering
			\includegraphics[height=6cm]{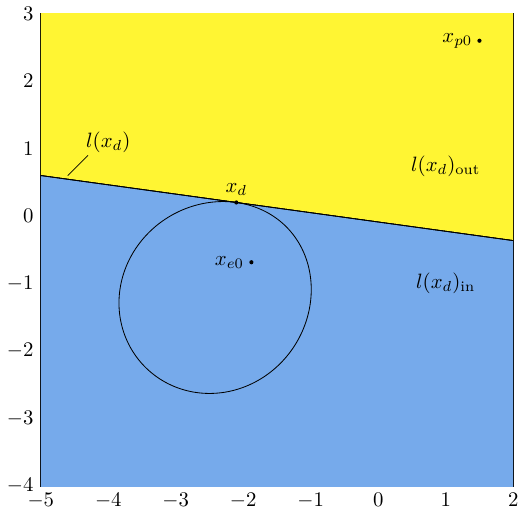}
			\caption{Take a point $x_d$ from $\partial \mathcal{D}_l (x_{p0},x_{e0})$. Draw the tangent $l(x_d)$ of $\partial \mathcal{D}_l (x_{p0},x_{e0})$ through $x_d$. According to the convexity of $\mathcal{D}_l (x_{p0},x_{e0})$, the whole of $\mathcal{D}_l (x_{p0},x_{e0})$ is located on the same side of $l(x_d)$.}
			\label{ztest0409_4} 
		\end{figure}\\
		Pick a point $x_d$ from $\partial \mathcal{D}_l (x_{p0},x_{e0})$. Draw the tangent $l(x_d)$ of $\partial \mathcal{D}_l (x_{p0},x_{e0})$ through $x_d$, as shown in Fig.~\ref{ztest0409_4}. According to the convexity of $\mathcal{D}_l (x_{p0},x_{e0})$, the whole of $\mathcal{D}_l (x_{p0},x_{e0})$ is located on the same side of $l(x_d)$. $\nabla_x\phi(x_d,0)$ is the normal direction of $l(x_d)$. The two sides of $l(x_d)$ can be represented as
		\begin{align}
			\left\lbrace x\in\mathbb{R}^2:\nabla_x\phi(x_d,0)\cdot (x-x_d)\le 0 \right\rbrace \label{tangent_out},\\
			\left\lbrace x\in\mathbb{R}^2:\nabla_x\phi(x_d,0)\cdot (x-x_d)\ge 0 \right\rbrace \label{tangent_in}.
		\end{align}
		Let \eqref{tangent_out} be the outside of $l(x_d)$, denoted as $l(x_d)_{\text{out}}$ (yellow region in Fig.~\ref{ztest0409_4}). Let \eqref{tangent_in} be the inside of $l(x_d)$, denoted as $l(x_d)_{\text{in}}$ (blue region in Fig.~\ref{ztest0409_4}). For any $t\in [t_0,t_f(\delta^*[u_e],u_e))$, consider 
		the optimization problem as follows:
		\begin{equation}\label{optimization problem}
			\begin{split}
				\min_{x\in\mathbb{R}^2} \nabla_x \phi (x_d,0)\cdot (x-x_d), \\
				\text{s.t.   } x\in \overline{\mathcal{D}_l(x_p(t;\delta^*[u_e]),x_e(t;u_e))}.
			\end{split}
		\end{equation}
		$ \nabla_x \phi (x_d,0)\cdot (x-x_d)$ is an affine function with respect to $x$.  $\overline{\mathcal{D}_l(x_p(t;\delta^*[u_e]),x_e(t;u_e))}$ is a bounded, closed, strictly convex set \cite{yan2021optimal}. The solution of the above optimization problem exists and is unique, which is denoted as $x_H(t)$. The optimal solution is also on $\partial\mathcal{D}_l(x_p(t;\delta^*[u_e]),x_e(t;u_e))$, i.e.
		\begin{align}\label{KKT_s}
			\phi(x_H(t),t)=0.
		\end{align} 
		By the KKT conditions, $x_H(t)$ satisfies
		\begin{align}
			\nabla_x\phi(x_d,0)-\lambda \nabla_x\phi(x_H(t),t)=0, \label{KKT_1}\\
			\lambda \phi(x_H(t),t)=0,\label{KKT_2}\\
			\lambda\ge 0.\label{KKT_3}
		\end{align}
		Take the derivative of equation \eqref{KKT_s} on both sides with respect to $t$:
		\begin{align}
			\frac{\partial\phi}{\partial t}(x_H(t),t)+\nabla_x\phi (x_H(t),t)\cdot \dot{x}_H(t)=0 \label{diff_t}.
		\end{align}
		Calculate $\frac{\partial\phi}{\partial t}(x_H(t),t)$,
		\begin{equation*}
			\begin{split}
				&\frac{\partial\phi}{\partial t}(x_H(t),t) \\
				=&\frac{x_p(t;\delta^*[u_e])-x_H(t)}{\left\| x_p(t;\delta^*[u_e])-x_H(t)\right\| }\cdot 	\alpha\delta^*[u_e](t)\\ &-\alpha \frac{x_e(t;u_e)-x_H(t)}{\left\| x_e(t;u_e)-x_H(t)\right\| }\cdot u_e(t)  \\
				=&-\alpha \left( \frac{x_H(t)-x_p(t;\delta^*[u_e])}{\left\| x_H(t)- x_p(t;\delta^*[u_e])\right\| }\cdot  \frac{x_c(t)-x_p(t;\delta^*[u_e])}{\left\| x_c(t)-x_p(t;\delta^*[u_e]) \right\| }\right. \\
				&\left. -\frac{x_H(t)-x_e(t;u_e)}{\left\|x_H(t)- x_e(t;u_e)\right\| }\cdot \frac{x_c(t)-x_e(t;u_e)}{\left\| x_c(t)-x_e(t;u_e) \right\| } \right),
			\end{split}
		\end{equation*}
		where $x_c(t)=x_c(x_p(t;\delta^*[u_e]),x_e(t;u_e),u_e(t))$. Then, by \thref{lemma1}, we have
		\begin{align}\label{co_lemma}
			\frac{\partial\phi}{\partial t}(x_H(t),t) \le 0.
		\end{align}
		Substituting \eqref{KKT_1}\eqref{KKT_3}\eqref{co_lemma} into \eqref{diff_t} and rearranging it, we obtain
		\begin{align*}
			\nabla_x\phi (x_d,0)\cdot \dot{x}_H(t)\ge 0.
		\end{align*}
		It is obvious that $x_H(0)=x_d$, then,
		\begin{align*}
			&\nabla_x\phi (x_d,0)\cdot (x_H(t)-x_d)\\
			=&\int_{0}^{t}\nabla_x\phi (x_d,0)\cdot \dot{x}_H(\tau) d\tau\ge 0.
		\end{align*}
		Thus, $\overline{\mathcal{D}_l(x_p(t;\delta^*[u_e]),x_e(t;u_e))}\subseteq l(x_d)_{\text{in}}$. Note that $x_d$ is picked arbitrarily. We obtain 
		\begin{align*}
			&\overline{\mathcal{D}_l(x_p(t;\delta^*[u_e]),x_e(t;u_e))} \\
			\subseteq &\bigcap_{x\in\partial \mathcal{D}_l (x_{p0},x_{e0})} l(x)_{\text{in}}=\overline{\mathcal{D}_l (x_{p0},x_{e0})}.
		\end{align*}
		where the last equality holds by \cite[Theorem 18.8]{rockafellar1997convex}. The proof is completed.\qed
	\end{pf}
	\begin{rem}
		When $X$ is a convex subset of $\mathbb{R}^2$, the conclusion in \thref{therom1} holds with just $\mathcal{D}_l(x_{p0},x_{e0}) \cap X$ replacing $\mathcal{D}_l(x_{p0},x_{e0})$ in \thref{therom1}(b).
	\end{rem}
	
	\section{Necessary condition for the existence of a PELIDR non-anticipative strategy with obstacles and zero capture radius}
	In the previous section, we proved that in the absence of obstacles, the pursuer has a PELIDR non-anticipative strategy. In next two sections, we consider this issue in the presence of obstacles. The obstacles considered are areas bounded by finite line segments, rays and lines. In this section, we present a counter example to illustrate that the pursuer's PELIDR non-anticipative strategies may not exist in the presence of obstacles. We also provide a necessary condition for the existence of a PELIDR non-anticipative strategy. Below are some basic relevant theories \cite{oyler2016pursuit}.
	\begin{thm}\thlabel{shortestline_form}
		In the presence of a set of polygonal obstacles, the shortest paths connecting any two points in $X$ are  broken lines, breaking at obstacle vertices.
	\end{thm}
	The proof of \thref{shortestline_form} can be found in \cite{chein1983routing}. 
	\begin{thm}\thlabel{domain_boundary_construction}
		In the presence of a set of polygonal obstacles, the boundary of the evader's dominance region consists of the following three types of curves:\\
		a) Cartesian ovals with focuses on the position of the pursuer or the evader or the vertexes of obstacles,\\
		b) Apollonius circles with focuses on the position of the pursuer or the evader or the vertexes of obstacles,\\
		c) circles centered on the vertexes of obstacles.
	\end{thm}
	\thref{domain_boundary_construction} is from \cite[Theorem 7]{oyler2016pursuit} and the discussion after its proof. \\ 
	The above conclusions still hold in the case that the obstacles are areas separated by finite line segments, rays, and lines. The proofs are the same as those given above.
	\begin{prop}\thlabel{partial_dL_norm}
		If $d_L$ is continuously differentiable near $(x_1,x_2)$, then 
		\begin{align*}
			\left\| \partial_{1}d_L(x_1,x_2)\right\|=\left\|  \partial_{2}d_L(x_1,x_2)\right\| =1,
		\end{align*}
		where $\partial_{i}d_L$ is the derivative of $d_L(x_1,x_2)$ with respect to $x_i$, $i=1,2$.
	\end{prop}
	The proof is presented in Appendix \ref{proof_partial_dL_norm}.

	\subsection{Counter example}
	In the presence of obstacles, the pursuer's PELIDR non-anticipative strategy does not always exist. We provide a counter example to illustrate this.
	\begin{exmp}\thlabel{con_example}
		\begin{figure}[h]
			\centering
			\includegraphics[height=6cm]{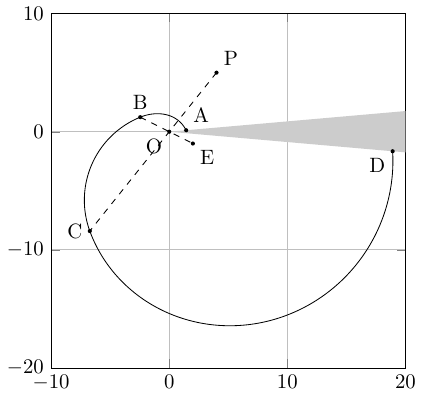} 
			\caption{Boundary of the evader's dominance region under the condition that $x_{p0}=(4,5)$, $x_{e0}=(2,-1)$, $\alpha=1.5$.}
			\label{figure_example_1}
		\end{figure}
		See Fig.~\ref{figure_example_1}. The boundary of the obstacle is an angle whose vertex is at the origin $O$. Each edge makes an angle of $5^{\circ}$ with the x-axis. The coordinate of the pursuer's initial position $P$ is $x_{p0}=(4,5)$. The coordinate of the evader's initial position $E$ is $x_{e0}=(2,-1)$. The speed ratio of pursuer to evader $\alpha$ is $1.5$. The boundary of the initial evader's dominance region $\partial \mathcal{D}(x_{p0},x_{e0})$ is shown in Fig.~\ref{figure_example_1}. $A$ and $D$ are intersections of the boundary of the initial evader's dominance region and the two edges of the obstacle. $B$ is the intersection of $\partial \mathcal{D}(x_{p0},x_{e0})$ and the extension line of $EO$. $C$ is the intersection of $\partial \mathcal{D}(x_{p0},x_{e0})$ and the extension line of $PO$. According to \thref{shortestline_form} and \thref{domain_boundary_construction}, $\partial \mathcal{D}(x_{p0},x_{e0})$ is divided into three part: the curve $AB$, the curve $BC$, the curve $CD$. The curve $AB$ is part of a Cartesian oval. The curve $BC$ is part of an Apollonius circle. The curve $CD$ is part of another Cartesian oval.\\
		\begin{figure}[h]
			\centering
			\includegraphics[height=6cm]{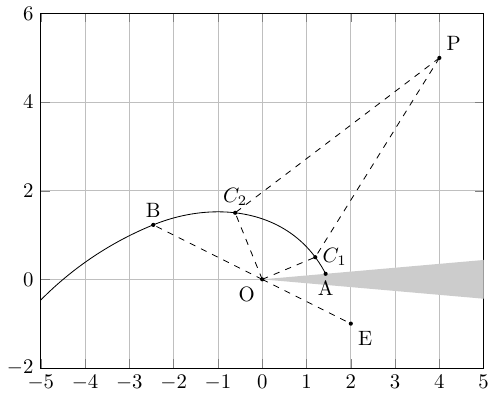} 
			\caption{Boundary of the evader's dominance region under the condition that $x_{p0}=(4,5)$, $x_{e0}=(2,-1)$, $\alpha=1.5$.}
			\label{figure_example_2}
		\end{figure}
		Let us consider the curve $AB$. See Fig.~\ref{figure_example_2}. Let $C_1,C_2$ be two different points on the curve $AB$. The shortest path connecting $E$ and $C_i$ is the broken line $EOC_i$, $i=1,2$. The shortest path connecting $P$ and $C_i$ is the line segment $PC_i$, $i=1,2$. Assume that there exists a PELIDR non-anticipative strategy for the pursuer. If the evader moves along the broken line $EOC_1$ by control input function $u_e$, the pursuer must move to $C_1$ directly to prevent the evader from entering the complement of $ \mathcal{D}(x_{p0},x_{e0})$. The pursuer's velocity direction is $\overrightarrow{PC_1}$ during a sufficiently small time interval $[t_0,t_0+\Delta t)$. Similarly, if the evader moves along the broken line $EOC_2$ by control input function $u_e'$, the pursuer must move to $C_2$ directly. The pursuer's velocity direction is  $\overrightarrow{PC_2}$ during the small time interval $[t_0,t_0+\Delta t)$. During the small time interval $[t_0,t_0+\Delta t)$, the evader's velocity direction is $\overrightarrow{EO}$ whether it is along broken line $EOC_1$ or broken line $EOC_2$, i.e.$u_e(t)=u_e'(t),t\in[t_0,t_0+\Delta t)$. However, the velocity directions of the pursuer are different,$\delta_p[u_e](t) \neq \delta_p[u_e'](t),t\in[t_0,t_0+\Delta t)$. It does not satisfy being non-anticipative \eqref{nonanticipate}. Thus a PELIDR non-anticipative strategy does not exist.  	
	\end{exmp}
	\subsection{Necessary condition}\label{subsec_Necessary_condition}
	Following the proof of \thref{therom1}, we analyze the evolution of the boundary of the evader's dominance region. Assume that the boundary of the evader's dominance region is smooth enough in a time interval $[t_1,t_2]$. Let
	\begin{align*}
		\phi(x,t)\triangleq  d_L(x,x_p(t))-\alpha d_L(x,x_e(t)). 
	\end{align*}
	The derivative of $\phi(x,t)$ with respect to $t$ is
	\begin{equation}\label{dphi_dt}
		\begin{split}
			&\frac{\partial \phi}{\partial t}(x,t)\\
			=& \partial_{2}d_L(x,x_p(t))\cdot\dot{x}_p(t)-\alpha \partial_{2}d_L(x,x_e(t))\cdot\dot{x}_e(t).
		\end{split}
	\end{equation}
	Assume that $\partial\mathcal{D}(x_p(t),x_e(t))$ has a smooth parametric representation $C(s,t)$, $C:I \times [t_1,t_2] \rightarrow \mathbb{R}^2$. We obtain
	\begin{align*}
		\phi(C(s,t),t)=0.
	\end{align*}
	Taking the derivative with respect to $t$ on both sides of the above equation, we obtain
	\begin{align}
		\nabla_x \phi(C(s,t),t) \cdot \frac{\partial C}{\partial t}(s,t)+\frac{\partial \phi}{\partial t}(C(s,t),t)=0	\label{total_differential}.
	\end{align}
	Combining \eqref{dphi_dt} and \eqref{total_differential}, 
	\begin{equation}\label{geometry_method_key}
		\begin{split}
			&\nabla_x \phi(C(s,t),t) \cdot \frac{\partial C}{\partial t}(s,t)\\ 
			=&\alpha \partial_{2}d_L(C(s,t),x_e(t))\cdot\dot{x}_e(t)\\
			&- \partial_{2}d_L(C(s,t),x_p(t))\cdot\dot{x}_p(t).	
		\end{split}
	\end{equation}
	Here, we explain the geometric meaning of \eqref{geometry_method_key}.\\ 
	Consider the left side of \eqref{geometry_method_key}. $\nabla_x \phi(C(s,t),t)$ represents the inner normal vector of the boundary curve  of the evader's dominance region. $\frac{\partial C}{\partial t}$ is the velocity of a point on the boundary. $\nabla_x \phi(C(s,t),t) \cdot \frac{\partial C}{\partial t}(s,t)>0$ means that the point represented by $s$ will move towards the inside of the evader's dominance region. $\nabla_x \phi(C(s,t),t) \cdot \frac{\partial C}{\partial t}(s,t)<0$ means that the point represented by $s$ will move towards the outside of the evader's dominance region. If $\nabla_x \phi(C(s,t),t) \cdot \frac{\partial C}{\partial t}(s,t)<0$ for all $s$, the boundary will expand outward.\\
	Consider the right side of \eqref{geometry_method_key}. $\partial_{1}d_L(x_1,x_2)$ and $ \partial_{2}d_L(x_1,x_2)$ represent the tangent vectors of the shortest path connecting $x_1,x_2$, at $x_1$ and $x_2$ respectively. The right side of \eqref{geometry_method_key} is related to the projections of the pursuer's and the evader's velocity vectors on the shortest paths connecting them and $C(s,t)$, respectively.\\
	\eqref{geometry_method_key} establishes the relationship between the movement of the evader's dominance region boundary and the movements of the pursuer and the evader. If a PELIDR non-anticipative strategy exists, the evader's dominance region cannot expand outward at the initial moment. From what we have discussed above, we provide the following necessary condition for the existence of a PELIDR non-anticipative strategy. 
	\begin{thm}\thlabel{necessary_condition}
		Given initial condition $t_0,\mathbf{x}_0$, the points $x_{c1},x_{c2}$ are on $ \partial \mathcal{D}(x_{p0},x_{e0})$. $d_L$ is continuously differentiable near $(x_{c1},x_{p0})$ and $(x_{c2},x_{p0})$. If $\exists \delta_p \in \Delta_p,\forall u_e \in \mathcal{U}_{t_0}^e(x_{e0}),\forall t \in [t_0,t_f(\delta_p[u_e],u_e)),x_e(t;t_0,x_{e0},u_e) \in {\mathcal{D}(x_{p0},x_{e0})}$, 
		then 
		\begin{equation}
			\begin{split}
				& \partial_{2}d_L(x_{c1},x_{p0})\cdot \partial_{2}d_L(x_{c2},x_{p0})\\
				&-  \partial_{2}d_L(x_{c1},x_{e0})\cdot \partial_{2}d_L(x_{c2},x_{e0})\ge0.
			\end{split}
		\end{equation}
	\end{thm}
	\begin{pf}
		For the sake of
		contradiction, assume that there exist $x_{c1},x_{c2} \in \partial \mathcal{D}(x_{p0},x_{e0})$ such that
		\begin{equation}\label{assumption}
			\begin{split}
				& \partial_{2}d_L(x_{c1},x_{p0})\cdot \partial_{2}d_L(x_{c2},x_{p0})\\
				-&  \partial_{2}d_L(x_{c1},x_{e0})\cdot \partial_{2}d_L(x_{c2},x_{e0})<0.
			\end{split}
		\end{equation}
		Select $u_e \in \mathcal{U}_e^{t_0}(x_{e0})$ that satisfies $u_e(t_0)= -\partial_{2}d_L(x_{c1},x_{e0})$. If $\delta_p[u_e](t_0) \neq  -\partial_{2}d_L(x_{c1},x_{p0})$, then
		\begin{equation}\label{ineq1}
			\begin{split}
				&\frac{\partial \phi}{\partial t}(x_{c1},t_0)\\
				=&\alpha\left(  \partial_{2}d_L(x_{c1},x_{p0}) \cdot \delta_p[u_e](t_0) \right. \\
				&-\left.  \partial_{2}d_L(x_{c1},x_{e0}) \cdot u_e(t_0) \right) \\
				=&\alpha\left(  \partial_{2}d_L(x_{c1},x_{p0}) \cdot \delta_p[u_e](t_0)  
				+ 1 \right)>0.
			\end{split}
		\end{equation}
		If $\delta_p[u_e](t_0) =  -\partial_{2}d_L(x_{c1},x_{p0})$, according to \eqref{assumption},
		\begin{equation}\label{ineq2}
			\begin{split}
				\frac{\partial \phi}{\partial t}(x_{c2},t_0)=&\alpha\left(  \partial_{2}d_L(x_{c2},x_{p0}) \cdot \delta_p[u_e](t_0) \right. \\
				-&\left.  \partial_{2}d_L(x_{c2},x_{e0}) \cdot u_e(t_0) \right) >0.
			\end{split}
		\end{equation}
		At least one of \eqref{ineq1} and \eqref{ineq2} is true. Without loss of generality, assume that \eqref{ineq1} is true. Take the Taylor expansion of $\phi$ at $(x_{c1},t_0)$ with respect to $t$:
		\begin{equation*}
			\begin{split}
				\phi(x_{c1},t_0+\Delta t)=\phi(x_{c1},t_0)+\frac{\partial \phi}{\partial t}(x_{c1},t_0)\Delta t+ o(\Delta t).
			\end{split}
		\end{equation*}
		Then, $\exists \tau>0,\forall \Delta t \in (0,\tau), \phi(x_{c1},t_0+\Delta t)>0$, i.e. $x_{c1} \in \mathcal{D}(x_p(t_0+\Delta t),x_e(t_0+\Delta t))$. Thus, $\exists x\in \mathcal{D}(x_{p0},x_{e0})^c, x\in \mathcal{D}(x_p(t_0+\Delta t),x_e(t_0+\Delta t)) $. According to \thref{reachable_is_domainance2}, the evader can reach $x$ before being captured with initial condition $t_0+\Delta t$ and $x_p(t_0+\Delta t),x_e(t_0+\Delta t)$, which is a contradiction. \qed
	\end{pf}
	\section{Sufficient condition for the existence of a PELIDR non-anticipative strategy with a single corner obstacle and zero capture radius}\label{a_simple_case}
	\begin{figure}[h]
		\centering
		\includegraphics[height=4cm]{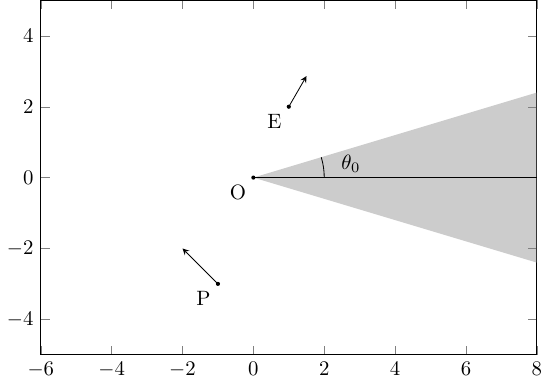}
		\caption{Pursuit-evasion games with a single corner obstacle}
		\label{obstacle x^+axis} 
	\end{figure}
	In this section, we focus on the case with a single corner obstacle. The obstacle is the gray area in Fig.~\ref{obstacle x^+axis}. Its boundary consists of two rays. The initial point of each ray is at the origin. The angle between the $x$-axis and each ray is $\theta_0$, $\theta_0\in(0,\pi/2)$. The region $X$ where the game takes place is as follows:
	\begin{align*}
		X=\left\lbrace (\rho\cos\theta,\rho\sin\theta):\rho\ge 0,\theta\in[\theta_0,2\pi-\theta_0]\right\rbrace. 
	\end{align*}
	This is similar with the ``configuration space" in \cite{zou2018optimal}. To the best of our knowledge, the case with a single corner obstacle is the simplest case with obstacles. We believe that the general cases can be locally approximated as the case with a single corner obstacle. Therefore, we think it is the basis for studying the general cases with obstacle. The details of the connection between the general cases and the case with a single corner obstacle is left in our future work. The following theorem is the core conclusion of this section.
	\begin{thm}\thlabel{simple_key_theorem}
		In the above circumstance, if the vertex of the obstacle is outside the closure of the evader's initial dominance region, i.e. $\left\| x_{p0}\right\| <\alpha\left\| x_{e0}\right\| $, then the pursuer has a PELIDR  non-anticipative strategy.
	\end{thm}
	
	\begin{rem}
		That the vertex of the obstacle is outside the closure of the evader's initial dominance region also implies that the pursuer can reach the vertex of the obstacle before the evader reaches it.
	\end{rem}
	A more rigorous statement of this theorem and its proof are given in Subsection \ref{proof sim key thm}. We first require some preliminaries to prove the theorem.
	\subsection{The shortest path distance}
	In this circumstance, the origin is the unique vertex of the obstacle. According to \thref{shortestline_form}, there are only two cases for the shortest paths connecting $x_1,x_2$ on $X$. \\
	(a)If the line segment whose endpoints are $x_1,x_2$ is contained in $X$, the shortest path connecting $x_1,x_2$ is the line segment.\\
	(b)If the line segment whose endpoints are $x_1,x_2$ is not contained in $X$, the shortest path connecting $x_1,x_2$ is the broken line connecting $x_1$, the origin, and $x_2$ in sequence.\\
	Let $(\rho_1,\theta_1),(\rho_2,\theta_2)$ denote the polar coordinates of $x_2,x_2$, respectively, where $\rho_i \ge 0 ,\theta_i \in [\theta_0,2\pi-\theta_0],i=1,2$. $\left| \theta_1-\theta_2 \right| > \pi$ is equivalent to that the line segment whose endpoints are $x_1,x_2$ intersects the obstacle. Thus, $d_L$ can be represented as follows:
	\begin{equation}\label{distance_of_shortest_path}
		d_L(x_1,x_2)=
		\begin{cases}
			\left\| x_1-x_2 \right\| &\left| \theta_1-\theta_2 \right|\le \pi\\
			\left\| x_1 \right\| + \left\| x_2 \right\| &\left| \theta_1-\theta_2 \right| > \pi   
		\end{cases}.
	\end{equation}
	The visible region of $x_1$ is
	\begin{align*}
		\mathcal{V}(x_1)=&\left\lbrace (\rho \cos \theta,\rho \sin \theta)\in X:
		\left| \theta-\theta_1\right| \le \pi \right\rbrace. 
	\end{align*}
	If $x_1$ is on the upper half plane, the visible region of $x_1$ can be represented as
	\begin{align*}
		\mathcal{V}(x_1)
		=\left\lbrace (\rho \cos \theta,\rho \sin \theta)\in X: \theta\in(\theta_0,\theta_1+\pi) \right\rbrace. 
	\end{align*}

	\subsection{Evader's dominance region and its boundary}
	Consider the evader's dominance region $\mathcal{D}(x_p,x_e) $ and its boundary $\partial \mathcal{D}(x_p,x_e)$ where $\left\| x_{p}\right\| <\alpha\left\| x_{e}\right\| $. The boundary is 
	\begin{align*}
		\partial \mathcal{D}(x_p,x_e) = \left\lbrace x\in X : d_L(x,x_p)-\alpha d_L(x,x_e) = 0  \right\rbrace.
	\end{align*}
	According to \thref{shortestline_in_domainance}, the shortest path connecting $x_e$ and any point on $\partial \mathcal{D}(x_p,x_e)$ does not pass the vertex of the obstacle when $\left\| x_{p}\right\| <\alpha\left\| x_{e}\right\| $. Thus, the shortest path is a line segment. Then, any point $x$ on $\partial \mathcal{D}(x_p,x_e)$ satisfies one of the following conditions:\\
	(a)The line segment whose endpoints are $x,x_p$ is contained in $X$. $x$ satisfies
	\begin{align}\label{A-circle}
		\left\| x-x_p \right\| - \alpha \left\| x-x_e\right\| =0.
	\end{align}
	\eqref{A-circle} represents the Apollonius circle whose focal points are $x_p,x_e$.\\
	(b)The line segment whose endpoints are $x,x_p$ is not contained in $X$. $x$ satisfies
	\begin{align}\label{C-oval}
		\left\| x \right\| - \alpha \left\| x-x_e\right\| = -\left\| x_p \right\|. 
	\end{align}
	\eqref{C-oval} represents the Cartesian oval of the second type whose focal points are the origin and $x_e$. In particular, $x\in\partial \mathcal{D}(x_p,x_e)$ satisfies both \eqref{A-circle} and \eqref{C-oval} if and only if $x,x_p$ and the origin are collinear. Let
	\begin{align*}
		\mathcal{A}(x_1,x_2) \triangleq \left\lbrace x\in \mathbb{R}^2 : \left\| x-x_1\right\|  - \alpha \left\| x-x_2\right\| > 0 \right\rbrace,\\
		\mathcal{C}(x_1,x_2) \triangleq \left\lbrace x\in \mathbb{R}^2 : \left\| x \right\| - \alpha \left\| x-x_2\right\| > -\left\| x_1 \right\| \right\rbrace.
	\end{align*}
	From what has been discussed above, we obtain 
	\begin{align}\label{boundary=A+D}
		\partial\mathcal{D}(x_p,x_e)\subseteq\partial\mathcal{A}(x_p,x_e)\cup\partial\mathcal{C}(x_p,x_e).
	\end{align}
	In fact, when $\left\| x_{p}\right\| <\alpha\left\| x_{e}\right\| $, we can define a set $\mathcal{F}(x_p,x_e)\subseteq \mathbb{R}^2 $ determined by $x_p,x_e$, whose boundary curve $\partial\mathcal{F}(x_p,x_e)$ satisfies the following properties:\\
	(a) $\partial\mathcal{F}(x_p,x_e) \subseteq \partial\mathcal{A}(x_p,x_e)\cup\partial\mathcal{C}(x_p,x_e)$;\\
	(b) $\partial\mathcal{D}(x_p,x_e)\subseteq\partial\mathcal{F}(x_p,x_e)$;\\
	(c) For any $u_e \in \partial B_2(\boldsymbol{0},1)$, draw a ray from the point $x_e$ along the vector $u_e$. There is a unique intersection point of the ray and $\partial \mathcal{F}(x_p,x_e)$;\\
	(d) $\partial\mathcal{F}(x_p,x_e)$ is a bounded, continuously differentiable, simple, closed curve.\\
	The definitions and properties of $\mathcal{F}(x_p,x_e)$ and $\partial\mathcal{F}(x_p,x_e)$ are presented in Appendix \ref{F_definition}. In detail, the definition of $\mathcal{F}(x_p,x_e)$ is from \eqref{eq_def_F} and \eqref{set_def_F}. (a) is from \eqref{eq_def_F} and \eqref{set_def_partialF}. (b) corresponds to \thref{D_in_F}. (c) corresponds to \thref{F-exist-uniq}. (d) corresponds to \thref{F_closed_curve}. 
	
	\subsection{Non-anticipative strategy $\delta^*$ and auxiliary strategy $\delta^\epsilon$}
	Next, we use $\partial \mathcal{F}(x_p,x_e)$ to construct the PELIDR non-anticipative strategy of the pursuer.\\
	Given $x_p,x_e \in X$ satisfying $\left\| x_p\right\| <\alpha\left\| x_e\right\|$, for any $u_e \in \partial B_2(\boldsymbol{0},1)$, take a ray from the point $x_e$ along the vector $u_e$. The ray and $\partial \mathcal{F}(x_p,x_e)$ have a unique intersection point, which is denoted by $x_c=x_c(x_p,x_e,u_e)$. Let $\gamma^*(x_p,x_e,u_e)$ be a continuous function that satisfies
	\begin{align*}
		\gamma^*(x_p,x_e,u_e)=\left\lbrace 
		\begin{aligned}
			&\frac{x_c-x_p}{\left\| x_c-x_p\right\| } 
			&x_c\in\partial\mathcal{A}(x_p,x_e)\\
			&-\frac{x_p}{\left\| x_p\right\| } 
			&x_c\in\partial\mathcal{C}(x_p,x_e)   
		\end{aligned}
		\right. ,
	\end{align*}
	when $\left\| x_p\right\| <\alpha\left\| x_e\right\|$.
	Given initial time $t_0$ and initial state $\mathbf{x}_0$, $\gamma^*$ can generate a non-anticipative strategy of the pursuer like \eqref{nonanti_stra2}. The non-anticipative strategy is denoted by $\delta_{t_0,\mathbf{x}_0}^*$.\\
	Because we can not directly prove that $\delta_{t_0,\mathbf{x}_0}^*$ is a PELIDR non-anticipative strategy, we introduce a group of auxiliary strategies.
	For $\epsilon>0$, let
	\begin{align*}
		\gamma^\epsilon(x_p,x_e,u_e)=\gamma^*(x_p,x_e,u_e)-\epsilon \partial_{1}d_L(x_p,x_e).
	\end{align*}
	Given initial time $t_0$ and initial state $\mathbf{x}_0$, for any $\epsilon$, $\gamma^\epsilon$ can generate a map from $\mathcal{U}_{t_0}^e$ to 
	\begin{align*}
		\left\lbrace u:\mathbb[t_0,+\infty)\rightarrow B_2(\boldsymbol{0},1+\epsilon):u \text{ is piecewise continuous} \right\rbrace 
	\end{align*}
	like \eqref{nonanti_stra2}, which is denoted by $\delta_{t_0,\mathbf{x}_0}^\epsilon$. This is non-anticipative \eqref{nonanticipate} and can be seen as a strategy of the pursuer. We will demonstrate that $\delta_{t_0,\mathbf{x}_0}^\epsilon$ is a PELIDR strategy for any $\epsilon>0$ in \thref{strictly}. Subsequently, we approximate $\delta_{t_0,\mathbf{x}_0}^*$ by $\delta_{t_0,\mathbf{x}_0}^\epsilon$ with $\epsilon\rightarrow0$ to prove that $\delta_{t_0,\mathbf{x}_0}^*$ is a PELIDR strategy in the proof of \thref{simple_key_theorem2}.\\
	We now present some properties of $\delta_{t_0,\mathbf{x}_0}^*$ and $\delta_{t_0,\mathbf{x}_0}^\epsilon$.
	\begin{lem}\thlabel{theo_of_dom}
		When $\left\| x_p\right\| <\alpha\left\| x_e\right\|$, $\forall x\in\partial \mathcal{D}(x_p,x_e)$, $\forall u_e\in\partial B_2(\boldsymbol{0},1)$, 
		\begin{align}\label{pro_of_dom}
			\partial_{2}d_L(x,x_e)\cdot u_e- \partial_{2}d_L(x,x_p)\cdot \gamma^*(x_p,x_e,u_e) \ge 0.
		\end{align}
	\end{lem}
	\begin{pf}
		It is easy to obtain that
		\begin{align*}
			& \partial_{2}d_L(x,x_e)=\frac{x_e-x}{\left\| x_e-x\right\| },\\
			& \partial_{2}d_L(x,x_p)=-\gamma^*\left(x_p,x_e,\frac{x-x_e}{\left\| x-x_e\right\| }\right),\\
			&u_e=\frac{x_c(x_p,x_e,u_e)-x_e}{\left\| x_c(x_p,x_e,u_e)-x_e\right\| }.
		\end{align*}
		Then, \eqref{pro_of_dom} converts to
		\begin{align*}
			&\gamma^*\left(x_p,x_e,\frac{x-x_e}{\left\| x-x_e\right\| }\right)
			\cdot\gamma^*\left(x_p,x_e,\frac{x_c-x_e}{\left\| x_c-x_e\right\| } \right)\\
			-&\frac{x-x_e}{\left\| x-x_e\right\| }\cdot\frac{x_c-x_e}{\left\| x_c-x_e\right\| } \ge 0.
		\end{align*}	
		Thus it is only necessary to prove that $\forall x_1,x_2\in \partial \mathcal{F}(x_p,x_e)$,
		\begin{equation}\label{cosine_eq}
			\begin{split}
				&\gamma^*\left(x_p,x_e,\frac{x_1-x_e}{\left\| x_1-x_e\right\| }\right)
				\cdot\gamma^*\left(x_p,x_e,\frac{x_2-x_e}{\left\| x_2-x_e\right\| } \right)\\
				-&\frac{x_1-x_e}{\left\| x_1-x_e\right\| }\cdot\frac{x_2-x_e}{\left\| x_2-x_e\right\| } \ge 0.
			\end{split}
		\end{equation}
		For any $x\in \partial \mathcal{F}(x_p,x_e)$, let $\chi=\chi(x)$ denote the angle between vectors $x-x_e$ and $x_p-x_e$ ($x_p-x_e$ is the starting edge, and $x-x_e$ is the ending edge, and counterclockwise is the positive direction). Let $\psi=\psi(x)$ denote the angle between vectors $\gamma^*\left(x_p,x_e,\frac{x-x_e}{\left\| x-x_e\right\| }\right)$ and $x_p-x_e$ ($x_p-x_e$ is the starting edge, $\gamma^*\left(x_p,x_e,\frac{x-x_e}{\left\| x-x_e\right\| }\right)$ is the ending edge, and counterclockwise is the positive direction). Because $\partial \mathcal{F}(x_p,x_e)$ is a simple closed curve, assume that $\chi(x_2)-\chi(x_1)\in[0,\pi]$ without loss of generality.
		Like in the proof of \thref{lemma1}, if
		\begin{align}\label{d_theta/d_eta}
			\left| \frac{ d\psi(x)}{d\chi(x)}\right| \le 1
		\end{align}
		holds on the whole closed curve $\partial \mathcal{F}(x_p,x_e)$,
		\begin{align*}
			&\left| \psi(x_2)-\psi(x_1)\right| =\left| \int_{\chi(x_1)}^{\chi(x_2)} \frac{ d\psi(x)}{d\chi(x)} d\chi(x)\right| \\
			\le& \int_{\chi(x_1)}^{\chi(x_2)} \left| \frac{ d\psi(x)}{d\chi(x)} \right|d\chi(x)  \le \int_{\chi(x_1)}^{\chi(x_2)} d\chi(x) \\
			\le& \chi(x_2)-\chi(x_1).
		\end{align*}
		Thus, 
		\begin{align*}
			\cos(\psi(x_2)-\psi(x_1)) \ge \cos(\chi(x_2)-\chi(x_1)).
		\end{align*}
		\eqref{cosine_eq} holds by the geometric interpretation of vectors' inner product. Let us prove that \eqref{d_theta/d_eta} holds on the whole of $\partial \mathcal{F}(x_p,x_e)$ when $\left\| x_p\right\| <\alpha\left\| x_e\right\|$.\\
		We have obtained that \eqref{d_theta/d_eta} holds on the whole Apollonius circle $\partial \mathcal{A}(x_p,x_e)$ in the proof of \thref{lemma1}. Thus, \eqref{d_theta/d_eta} is true if $x \in \partial\mathcal{A}(x_p,x_e)$.\\
		When $x\in \partial\mathcal{C}(x_p,x_e)$,
		\begin{align*}
			\gamma^*\left(x_p,x_e,\frac{x-x_e}{\left\| x-x_e\right\| }\right)=-\frac{x_p}{\left\| x_p \right\|}.
		\end{align*}
		$\psi$ is invariable but $\chi$ is variable while $x$ moves on $\partial\mathcal{C}(x_p,x_e)$. Thus
		\begin{align*}
			\left| \frac{ d\psi(x)}{d\chi(x)}\right| = 0.
		\end{align*}
		Therefore, \eqref{d_theta/d_eta} holds on the whole of $\partial \mathcal{F}(x_p,x_e)$. \qed
	\end{pf}
	\begin{prop}\thlabel{lemma_of_increament}
		When $\left\| x_p\right\| <\alpha\left\| x_e\right\|$, $\forall x \in \partial \mathcal{D}(x_p,x_e)$,
		\begin{align}\label{increment}
			\partial_{1}d_L(x_p,x_e) \cdot  \partial_{2}d_L(x,x_p)>0.
		\end{align}
	\end{prop}
	The proof is presented in Appendix \ref{proof_lemma_of_increament}.
	
	\subsection{Proof of \thref{simple_key_theorem}} \label{proof sim key thm}
	Before presenting the rigorous statement and the proof of \thref{simple_key_theorem}, we state a lemma.
	\begin{lem}\thlabel{strictly}
		In the circumstance described at the beginning of this section, if the initial state $\mathbf{x}_0=(x_{p0},x_{e0})$ satisfies $\left\| x_{p0}\right\| <\alpha\left\| x_{e0}\right\|$ and the pursuer uses the strategy $\delta_{t_0,\mathbf{x}_0}^\epsilon$($\epsilon>0$), then for any $u_e\in \mathcal{U}_e^{t_0}(x_{e0})$, for any $t\in[t_0,t_f(\delta_{t_0,\mathbf{x}_0}^\epsilon[u_e],u_e))$, $\mathcal{D}(x_p(t),x_e(t))\subseteq \mathcal{D}(x_{p0},x_{e0})$.
	\end{lem}
	\begin{pf}
		Pick $u_e$ from $\mathcal{U}_e^{t_0}$. In this proof, we abbreviate \\$x_p(t;t_0,x_{p0},\delta_{t_0,\mathbf{x}_0}^*[u_e])$ and $x_e(t;t_0,x_{e0},u_e)$ as $x_p(t)$ and $x_e(t)$, respectively. For the sake of contradiction, assume that $\exists u_e\in\mathcal{U}_e^{t_0}(x_{e0})$, $\exists t\in[t_0,t_f(\delta_{t_0,\mathbf{x}_0}^\epsilon[u_e],u_e))$, $\exists x\in \mathcal{D}(x_p(t),x_e(t))$, $x\notin \mathcal{D}(x_{p0},x_{e0})$. Let
		\begin{align*}
			\tau^* \triangleq \inf\left\lbrace t:\exists x \in \mathcal{D}(x_p(t),x_e(t)), x\notin \mathcal{D}(x_{p0},x_{e0}) \right\rbrace. 
		\end{align*}
		Next, we show that $\mathcal{D}(x_p(\tau^*),x_e(\tau^*))\subseteq \mathcal{D}(x_{p0},x_{e0})$. For the sake of contradiction, assume that there exists $x$ in $\mathcal{D}(x_p(\tau^*),x_e(\tau^*))$ such that $x \notin \mathcal{D}(x_{p0},x_{e0})$. By the definition,
		\begin{align*}
			d_L(x,x_p(\tau^*))-\alpha d_L(x,x_e(\tau^*))>0.
		\end{align*}
		According to the continuities of $d_L$, $x_p(\cdot)$ and $x_e(\cdot)$, there exists  $\Delta t>0$ such that 
		\begin{align*}
			d_L(x,x_p(\tau^*-\Delta t))-\alpha d_L(x,x_e(\tau^*-\Delta t))>0,
		\end{align*}
		i.e. $x\in\mathcal{D}(x_p(\tau^*-\Delta t),x_e(\tau^*-\Delta t))$. This contradicts the definition of $\tau^*$.
		Thus, we obtain $\mathcal{D}(x_p(\tau^*),x_e(\tau^*))\subseteq \mathcal{D}(x_{p0},x_{e0})$. According to the definition of infimum, there exist sequences $\left\lbrace x_n\right\rbrace _{n=1}^{+\infty} \subseteq X$ and $\left\lbrace t_n\right\rbrace _{n=1}^{+\infty} \subseteq (\tau^*,+\infty)$ such that
		\begin{align}
			&x_n\in\mathcal{D}(x_p(t_n),x_e(t_n)),\label{in_C_t}\\
			&x_n\notin \mathcal{D}(x_{p0},x_{e0}),\label{not_in_C_0}\\
			&\lim_{n\rightarrow +\infty}t_n=\tau^*.\label{lim_t}
		\end{align}
		$\left\lbrace x_n\right\rbrace _{n=1}^{+\infty} $ is bounded, thus it contains a convergent subsequence. Let the subsequence be itself without loss of generality. Let
		\begin{align}\label{lim_x}
			x^*\triangleq\lim_{n\rightarrow +\infty}x_n.
		\end{align}
		According to \eqref{in_C_t} and \eqref{not_in_C_0},
		\begin{align*}
			d_L(x_n,x_p(t_n))-\alpha d_L(x_n,x_e(t_n))>0,\\
			d_L(x_n,x_{p0})-\alpha d_L(x_n,x_{e0})\le 0.
		\end{align*}
		Taking $n\rightarrow +\infty $,
		\begin{align}
			d_L(x^*,x_p(\tau^*))-\alpha d_L(x^*,x_e(\tau^*))\ge 0,\label{x*_ge_0}\\
			d_L(x^*,x_{p0})-\alpha d_L(x^*,x_{e0})\le 0.\label{x*_le_0}
		\end{align}
		According to $\mathcal{D}(x_p(\tau^*),x_e(\tau^*))\subseteq \mathcal{D}(x_{p0},x_{e0})$, we obtain
		\begin{align}
			\overline{\mathcal{D}(x_p(\tau^*),x_e(\tau^*))} \subseteq \overline{\mathcal{D} (x_{p0},x_{e0})},\label{closure_subset}\\
			\mathcal{D}(x_{p0},x_{e0})^c \subseteq \mathcal{D}(x_p(\tau^*),x_e(\tau^*))^c.\label{complement_subset}
		\end{align}
		By\eqref{x*_ge_0} and \eqref{closure_subset}, we have
		\begin{align*}
			d_L(x^*,x_{p0})-\alpha d_L(x^*,x_{e0})\ge 0.
		\end{align*}
		Combining the above inequality and\eqref{x*_le_0},
		\begin{align*}
			d_L(x^*,x_{p0})-\alpha d_L(x^*,x_{e0})= 0.
		\end{align*}
		By \eqref{x*_le_0} and \eqref{complement_subset}, we obtain
		\begin{align*}
			d_L(x^*,x_p(\tau^*))-\alpha d_L(x^*,x_e(\tau^*))\le 0.
		\end{align*}
		Combining the above inequality and \eqref{x*_ge_0}, 
		\begin{align*}
			d_L(x^*,x_p(\tau^*))-\alpha d_L(x^*,x_e(\tau^*))= 0.
		\end{align*}
		Thus,
		\begin{align*}
			x^*\in \partial\mathcal{D}(x_{p0},x_{e0})\cap\partial\mathcal{D}(x_p(\tau^*),x_e(\tau^*)).
		\end{align*}
		Let
		\begin{align*}
			\phi (x,t)=d_L(x,x_p(t))-\alpha d_L(x,x_e(t)).
		\end{align*}
		The derivative of $\phi$ with respect to $t$ is
		\begin{align*}
			\frac{\partial \phi}{\partial t}(x,t)=&\alpha  \partial_{2}d_L(x,x_p(t)) \cdot \delta_{t_0,\mathbf{x}_0}^\epsilon[u_e](t)\\
			&-\alpha  \partial_{2}d_L(x,x_e(t)) \cdot u_e(t).
		\end{align*}
		Due to $x^*\in \partial\mathcal{D}(x_p(\tau^*),x_e(\tau^*))$, by \thref{theo_of_dom} and \thref{lemma_of_increament}, we obtain
		\begin{align*}
			&\frac{\partial \phi}{\partial t}(x^*,\tau^*)\\
			=&\alpha  \partial_{2}d_L(x^*,x_p(\tau^*)) \cdot \gamma(x_p(\tau^*),x_e(\tau^*),u_e(\tau^*))\\
			&-\epsilon \alpha  \partial_{2}d_L(x^*,x_p(\tau^*)) \cdot \partial_{1}d_L(x_p(\tau^*),x_e(\tau^*))\\
			&-\alpha  \partial_{2}d_L(x^*,x_e(\tau^*)) \cdot u_e(\tau^*)<0.
		\end{align*}
		Due to the continuity of $\frac{\partial \phi}{\partial t}$, $\exists M>0$, $\exists \delta_1>0$, $\forall x\in B_2(x^*,\delta_1)$, $\forall t \in (\tau^*,\tau^*+\delta_1)$, $\frac{\partial \phi}{\partial t}(x,t)<-M$.\\
		According to the Cauchy mean value theorem, $\exists \delta_2>0$, $\forall x \in B_2(x^*,\delta_2)$, $\forall t\in(\tau^*,\tau^*+\delta_2)$, $\exists \Theta_{x,t}\in(\tau^*,t)$,
		\begin{align*}
			\phi(x,t)-\phi(x,\tau^*)=\frac{\partial\phi}{\partial t}(x,\Theta_{x,t})(t-\tau^*).
		\end{align*}
		Let $\delta_3=\min\left\lbrace \delta_1,\delta_2\right\rbrace  $. When $x \in B_2(x^*,\delta_3)\cap\mathcal{D}(x_p(\tau^*),x_e(\tau^*))^c$, $t\in(\tau^*,\tau^*+\delta_3)$, we have
		\begin{align*}
			\phi(x,t)&=\phi(x,\tau^*)+\frac{\partial\phi}{\partial t}(x,\Theta_{x,t})(t-\tau^*)\\
			&<-M(t-\tau^*)<0
		\end{align*}
		i.e.
		\begin{align*}
			x\notin\mathcal{D}(x_p(t),x_e(t)).
		\end{align*}
		This contradicts \eqref{in_C_t}-\eqref{lim_x}. \qed
	\end{pf}
	We present a more detailed statement of \thref{simple_key_theorem} and prove it.
	\begin{thm}[Detailed statement of \thref{simple_key_theorem}]\thlabel{simple_key_theorem2}
		In the circumstance described at the beginning of this section, if the initial state $\mathbf{x}_0=(x_{p0},x_{e0})$ satisfies $\left\| x_{p0}\right\| <\alpha\left\| x_{e0}\right\|$ and the pursuer uses the strategy $\delta_{t_0,\mathbf{x}_0}^*$, then\\
		(a)for any $u_e\in \mathcal{U}_e^{t_0}(x_{e0})$ and $t\in[t_0,t_f(\delta_{t_0,\mathbf{x}_0}^*[u_e],u_e))$, \\$x_e(t;t_0,x_{e0},u_e)\in \mathcal{D}(x_{p0},x_{e0})$;\\
		(b)there exists a constant $c$, for any $u_e\in \mathcal{U}_e^{t_0}(x_{e0})$, $t_f(\delta_{t_0,\mathbf{x}_0}^*[u_e],u_e)-t_0<c$.
	\end{thm}
	\begin{pf}
		(a) Pick $u_e(\cdot)$ from $\mathcal{U}_e^{t_0}$. In this proof, we abbreviate $x_p(t;t_0,x_{p0},\delta_{t_0,\mathbf{x}_0}^*[u_e])$ and $x_e(t;t_0,x_{e0},u_e)$ as $x_p(t)$ and $x_e(t)$, respectively.
		It is only necessary to prove that $\mathcal{D}(x_p(t),x_e(t))\subseteq \mathcal{D}(x_{p0},x_{e0}) $ for any $u_e\in \mathcal{U}_e^{t_0}(x_{e0})$ and $t\in[t_0,t_f(\delta_{t_0,\mathbf{x}_0}^*[u_e],u_e))$.
		For the sake of contradiction, assume that there exists $u_e\in\mathcal{U}_e^{t_0}(x_{e0})$ such that $\exists t\in[t_0,t_f(\delta_{t_0,\mathbf{x}_0}[u_e],u_e))$, $\exists x\in \mathcal{D}(x_p(t),x_e(t))$, $x\notin \mathcal{D}(x_{p0},x_{e0})$. Let
		\begin{align*}
			\tau^* \triangleq \inf\left\lbrace t:\exists t \in \mathcal{D}(x_p(t),x_e(t)), x\notin \mathcal{D}(x_{p0},x_{e0}) \right\rbrace 
		\end{align*}
		Similarly to the proof of \thref{strictly}, we have\\ $\mathcal{D}(x_p(\tau^*),x_e(\tau^*))\subseteq \mathcal{D}(x_{p0},x_{e0})$ and there exist sequences $\left\lbrace x_n\right\rbrace _{n=1}^{+\infty} \subseteq X$ and $\left\lbrace t_n\right\rbrace _{n=1}^{+\infty} \subseteq (\tau^*,+\infty)$ such that
		\begin{align}
			&x_n\in\mathcal{D}(x_p(t_n),x_e(t_n)),\label{in_C_t_2}\\
			&x_n\notin \mathcal{D}(x_{p0},x_{e0}),\label{not_in_C_0_2}\\
			&\lim_{n\rightarrow +\infty}t_n=\tau^*,\label{lim_t_2}\\
			&\lim_{n\rightarrow +\infty}x_n=x^*.\label{lim_x_2}
		\end{align}
		Let the pursuer use the non-anticipative strategy $\delta_{t_0,\mathbf{x}_0}^*$ in time interval $[t_0,\tau^*)$ and use the non-anticipative strategy $\delta_{\tau^*,\mathbf{x}(\tau^*)}^\epsilon$ after $\tau^*$. For $t>\tau^*$, let  
		\begin{align*}
			x_p^\epsilon(t)=x_p(t;\tau^*,x_p(\tau^*),\delta_{\tau^*,\mathbf{x}(\tau^*)}^\epsilon[u_e]).
		\end{align*}
		According to the continuous dependence theory of ordinary differential equations \cite{hale2009ordinary}\cite{hsieh1999basic}, $\exists\Delta t>0, \forall t\in(\tau^*,\tau^*+\Delta t)$, 
		\begin{align*}
			\lim_{\epsilon\rightarrow 0}x_p^\epsilon(t)=x_p(t).
		\end{align*}
		Select $x\in \mathcal{D}(x_{p0},x_{e0})^c$. $x$ satisfies $x \in \mathcal{D}(x_p(\tau^*),x_e(\tau^*))^c$. By \thref{strictly}, $\forall t \in (\tau^*,\tau^*+\Delta t)$, $x\in\mathcal{D}(x_p^\epsilon(t),x_e(t))^c$, i.e.
		\begin{align*}
			d_L(x,x_p^\epsilon(t))-\alpha d_L(x,x_e(t))\le 0.
		\end{align*}
		Taking $\epsilon\rightarrow0$, 
		\begin{align*}
			d_L(x,x_p(t))-\alpha d_L(x,x_e(t))\le 0,
		\end{align*}
		i.e. $x\notin \mathcal{D}(x_p(t),x_e(t))$, which contradicts \eqref{in_C_t_2}-\eqref{lim_x_2}. \\
		(b) Taking the derivative of $d_L(x_p(t),x_e(t))$ with respect to $t$,
		\begin{align*}
			&\frac{d}{dt}d_L(x_p(t),x_e(t))\\
			=&\alpha \partial_{1}d_L(x_p(t),x_e(t))\cdot \gamma^*(x_p(t),x_e(t),u_e(t))\\
			&+ \partial_{2}d_L(x_p(t),x_e(t))\cdot u_e(t).
		\end{align*}
		By (a), we obtain that $\forall t\in [t_0,t_f)$, $\left\| x_p(t)\right\|<\alpha\left\| x_e(t)\right\|$. We will prove that $\frac{d}{dt}d_L(x_p(t),x_e(t)) \le 1-\alpha$ when $\left\| x_p(t)\right\|<\alpha\left\| x_e(t)\right\| $ through the following 4 cases.
		
		Case 1: $x_c(x_p(t),x_e(t),u_e(t))\in \partial\mathcal{A}(x_p,x_e)$ and $x_p,x_e$ are visible to each other. By the proof of \thref{therom1}(a),
		\begin{align*}
			\frac{d}{dt}d_L(x_p(t),x_e(t))\le 1-\alpha.
		\end{align*}
		\begin{figure}[h]
			\centering
			\includegraphics[height=6cm]{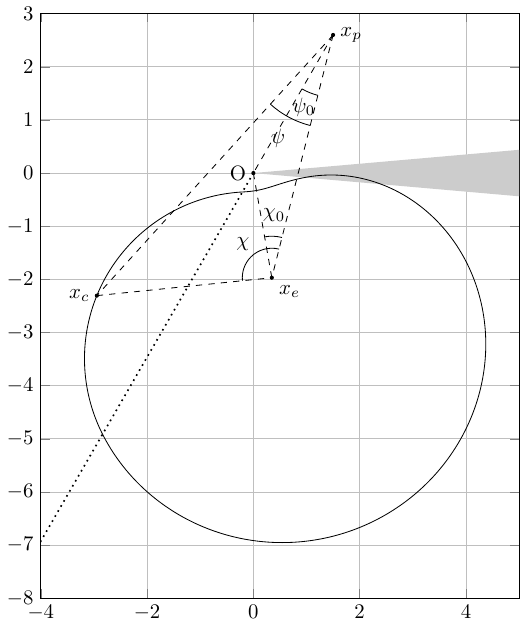}
			\caption{In the proof of \thref{simple_key_theorem2}, $x_p,x_e$ are not visible to each other and $x_c \in \partial \mathcal{A}(x_p,x_e)$.}
			\label{ztest0409_1}
		\end{figure}\\
		Case 2: $x_c(x_p(t),x_e(t),u_e(t))\in \partial\mathcal{A}(x_p,x_e)$ and $x_p,x_e$ are not visible to each other. As shown in Fig.~\ref{ztest0409_1}, let $\psi$ denote the angle between vectors $x_c-x_p$ and $x_e-x_p$, $\psi_0$ denote the angle between vectors $-x_p$ and $x_e-x_p$, $\chi$ denote the angle between vector $x_c-x_e$ and vector $x_p-x_e$, and $\chi_0$ denote the angle between vectors $-x_e$ and $x_p-x_e$. It is easy to obtain that
		\begin{align*}
			&\frac{d}{dt}d_L(x_p(t),x_e(t))\\
			=&-\alpha\cos(\psi-\psi_0)-\cos(\chi-\chi_0)\\
			\le&-\alpha\cos\psi-\cos\chi.
		\end{align*}
		By the proof of \thref{therom1}(a), $-\alpha\cos\psi-\cos\chi\le1-\alpha$.
		
		Case 3: $x_c(x_p(t),x_e(t),u_e(t))\in \partial\mathcal{C}(x_p,x_e)$ and $x_p,x_e$ are visible to each other. We obtain 
		\begin{align*}
			\partial_{1}d_L(x_p(t),x_e(t))=\frac{x_p-x_e}{\left\| x_p-x_e\right\| },\\
			\gamma^*(x_p(t),x_e(t),u_e(t))=-\frac{x_p}{\left\| x_p\right\|},\\
			\partial_{2}d_L(x_p(t),x_e(t))=\frac{x_e-x_p}{\left\| x_e-x_p\right\| },\\
			u_e(t)=\frac{x_c-x_e}{\left\| x_c-x_e\right\| }.\\
		\end{align*}
		\begin{figure}[h]
			\centering
			\includegraphics[height=6cm]{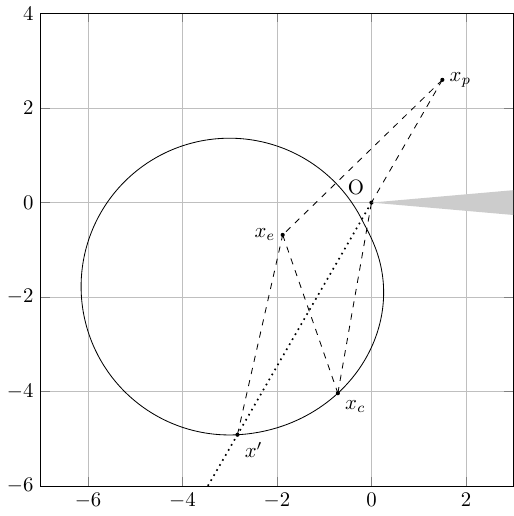}
			\caption{In the proof of \thref{simple_key_theorem2}, $x_p,x_e$ is visible to each other and $x_c \in \partial \mathcal{C}(x_p,x_e)$.}
			\label{ztest0409_2} 
		\end{figure}\\
		As shown in Fig.~\ref{ztest0409_2}, there exists two intersection points of $\partial\mathcal{F}(x_p,x_e)$ and the straight line passing through $x_p$ and the origin. The two intersection points are on $\partial\mathcal{A}(x_p,x_e)$. Let $x'$ denote the farthest intersection point from $x_p$. $x'$ satisfies
		\begin{align*}
			\frac{x'-x_p}{\left\| x'-x_p\right\| }=-\frac{x_p}{\left\| x_p\right\|},
		\end{align*}
		and
		\begin{align*}
			\frac{x_e-x_p}{\left\| x_e-x_p\right\| }\cdot\frac{x'-x_e}{\left\| x'-x_e\right\| }\ge\frac{x_e-x_p}{\left\| x_e-x_p\right\| }\cdot\frac{x_c-x_e}{\left\| x_c-x_e\right\| }.
		\end{align*}
		Combining the above and the proof of \thref{therom1}(a), we have
		\begin{align*}
			&\frac{d}{dt}d_L(x_p(t),x_e(t))\\
			=&\alpha \frac{x_p-x_e}{\left\| x_p-x_e\right\| }\cdot\frac{x'-x_p}{\left\| x'-x_p\right\| }
			+\frac{x_e-x_p}{\left\| x_e-x_p\right\| }\cdot\frac{x_c-x_e}{\left\| x_c-x_e\right\| }\\
			\le &\alpha \frac{x_p-x_e}{\left\| x_p-x_e\right\| }\cdot\frac{x'-x_p}{\left\| x'-x_p\right\| }
			+\frac{x_e-x_p}{\left\| x_e-x_p\right\| }\cdot\frac{x'-x_e}{\left\| x'-x_e\right\| }\\
			\le&1-\alpha.
		\end{align*}
		Case 4: $x_c(x_p(t),x_e(t),u_e(t))\in \partial\mathcal{C}(x_p,x_e)$ and $x_p,x_e$ are not visible to each other. We have
		\begin{align*}
			&\partial_{1}d_L(x_p(t),x_e(t))=- \partial_{2}d_L(x_p(t),x_e(t))\\
			=&-\gamma^*(x_p(t),x_e(t),u_e(t))=\frac{x_p}{\left\| x_p\right\| }.
		\end{align*}
		Thus,
		\begin{align*}\\
			&\frac{d}{dt}d_L(x_p(t),x_e(t))\\
			=&-\alpha + \partial_{2}d_L(x_p(t),x_e(t),)\cdot u_e(t)\le 1-\alpha.
		\end{align*}
		In summary, $\frac{d}{dt}d_L(x_p(t),x_e(t)) \le 1-\alpha$, when $\left\| x_p(t)\right\|<\alpha\left\| x_e(t)\right\|  $. Thus 
		\begin{align*}
			t_f(\delta^*[u_e],u_e)&\le t_0+\frac{\left\| x_{p0}-x_{e0}\right\|-l}{\alpha-1}.
		\end{align*}
		The proof is completed. \qed
	\end{pf}

	\section{Application}
	\begin{figure}[h]
		\centering
		\includegraphics[height=4cm]{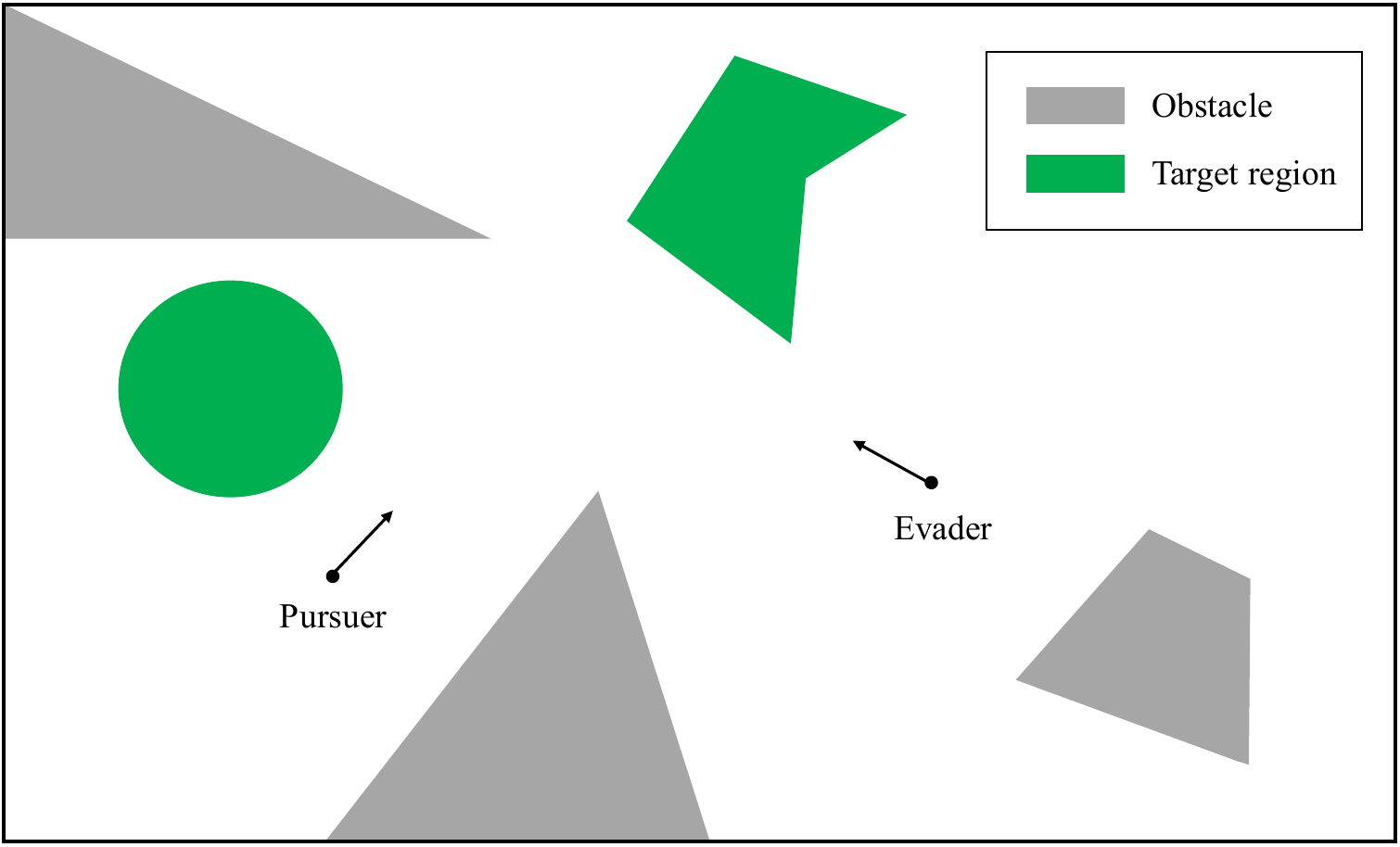}
		\caption{Target defence games.}
		\label{target_defence} 
	\end{figure}
	In this section, we discuss the application of the evader's dominance region in static target defense games.\\
	Based on the pursuit-evasion game described in Section \ref{PEgame}, a set $T$ is included in $X$, which is called target region, as shown in Fig.~\ref{target_defence}. The initial position of the evader is outside the target region. The problem is whether the pursuer has a non-anticipative strategy by which the pursuer can capture the evader before it reaches the target region $T$ regardless of the evader's control input function.\\
	In the absence of obstacles, this problem is perfectly solved. Such a strategy exists if and only if 
	\begin{align}\label{intesection_not_exists}
		\mathcal{D}(x_{p0},x_{e0})\cap T=\emptyset,
	\end{align}
	according to \thref{reachable_is_domainance1} and \thref{therom1}. Such a non-anticipative strategy is as defined in Section \ref{section without obstacle} if it exists.\\
	In the presence of obstacles, $\mathcal{D}_l(x_{p0},x_{e0})\cap T=\emptyset$ is a sufficient but not necessary condition for the existence of such a strategy. The sufficiency can be obtained from \thref{reachable_is_domainance2}. The non-necessity is explained as follows. Consider the case in \thref{con_example}. Let $T=\overline{\mathcal{D} (x_{p0},x_{e0})}^c$. Then \eqref{intesection_not_exists} holds. However, we have illustrated in \thref{con_example} that the pursuer does not have a non-anticipative strategy by which it can guarantee the prevention of the evader reaching the target region.
	
	\section{Conclusion}
	In this article, we reveal the properties of the evader's dominance region in pursuit-evasion games. We show that the evader's dominance region is its reachable region in the open-loop sense. We explore conditions under which the pursuer has a non-anticipative strategy to prevent the evader leaving its initial dominance region before capture. We prove that such a non-anticipative strategy exists in the absence of obstacles. We provide a counter example to illustrate that such a non-anticipative strategy does not always exist in the presence of obstacles and present a necessary condition for the existence of such a non-anticipative strategy. In the case of a single corner obstacle, we provide a sufficient condition for the existence of this non-anticipative strategy. We also discuss the application of the dominance region in the target defense game under the non-anticipative information pattern. This problem is perfectly solved in the absence of obstacles.\\ 
	We think the limitation of our work is that the sufficient condition is derived in an overly special case. General cases should be considered in future works.\\
	Compared to solving numerical solution of HJI equation, geometric method is effective for handle multi-agent game. The difficulties may be the task assignment for the agents when the scales is too large. In some special circumstances (such as convex target \cite{yan2022matching}), the efficiency of task assignment algorithms can be improved. Multi-agent games with obstacles deserve further study.\\
	Geometric methods are also possible to be applied in more complex dynamical systems. In \cite{yan2023homicidal}\cite{yan2024multiplayer}, the authors address static target defense games with dubins car defender. They proposed two defender's strategies inspired by the evader's dominance region. It is also worth studying how to transfer these strategies to cases with obstacles.

	\bibliographystyle{unsrt}        
	\bibliography{autosam}           
	
	
	
	\appendix
	
	\section{Challenge without the access to the opponent's current control input}\label{challenge}
	\begin{figure}[h]
		\centering
		\includegraphics[height=5cm]{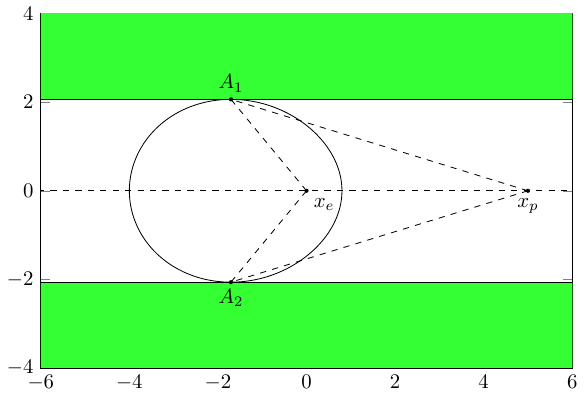}
		\caption{The boundary of the evader's dominance region is tangent to the boundary of the target at two points.}
		\label{zc_oval_revised1} 
	\end{figure}
	
	This example resembles to the problem discussed in \cite[Appendix A]{dorothy2024one}. See Fig.~\ref{zc_oval_revised1}. There is no obstacle on the plane. The target (green regions) is the union of two disjoint half planes. The boundary of each half plane is parallel to the symmetry axis of the evader's dominance region and tangent to the boundary of the dominance region. The tangent points are denoted by $A_1,A_2$. If the evader goes to $A_1$ and the pursuer does not know the evader's current velocity, the evader may choose to go to $A_2$. This causes that the evader passes through $A_1$ and enters the target according to the definition of dominance regions.  Therefore, the pursuer can not guarantee prevention of the evader entering the target if it does not know the evader's current velocity.
	
	\section{Proof of \thref{partial_dL_norm}}\label{proof_partial_dL_norm}
	\begin{pf}
		Let $s_1,s_2,...s_k$ denote the vertexes of obstacles. $s_0=x_2$. Let $S=\left\lbrace s_0, s_1,...,s_k\right\rbrace$. 
		Consider the shortest path map (partition) of $S$ \cite{lee1984euclidean}: 
		$R(s_0),R(s_1),...,R(s_k)$. 
		For any $i$, for any $x$ in the interior of $R(s_i)$, $d_L(x,x_2)=\left\| x-s_i\right\|+d_L(s_i,x_2)$. Then 
		\begin{align*}
			\left\| \partial_{1}d_L(x,x_2)\right\|=\left\| \frac{x_1-s_i}{\left\| x_1-s_i\right\| }\right\| =1,
		\end{align*}
		for any $x$ in the interior of $R(s_i)$. $\left\| \partial_{1}d_L(x,x_2)\right\|=1$ holds on the boundary of $R(s_i)$ by its continuous differentiability. We have completed the proof of $\left\| \partial_{1}d_L(x,x_2)\right\|=1$. $\left\|  \partial_{2}d_L(x_1,x_2)\right\|=1$ is obtained by the symmetry of $d_L$.\qed
	\end{pf}
	\section{Definitions and properties of $\mathcal{F}(x_p,x_e)$ and $\partial\mathcal{F}(x_p,x_e)$}\label{F_definition}
	In this section, we present the definitions and properties of $\mathcal{F}(x_p,x_e)$ and $\partial\mathcal{F}(x_p,x_e)$ in Section \ref{a_simple_case}. All $(x_p,x_e)$ in this section satisfy $\left\| x_p\right\|\le \alpha \left\| x_e \right\|  $. With $I$ denoting an interval or the empty set, let
	\begin{align*}
		\mathcal{G}(I)=\left\lbrace (\rho\cos\theta,\rho\sin\theta):\rho\ge0,\theta\in I\right\rbrace.
	\end{align*}
	\begin{prop}\thlabel{uniq_A-circle}
		For any $x\in\overline{\mathcal{A}(x_p,x_e)}$, select a point $x'$ from the line segment whose endpoints are $x,x_e$. Then, $x'\in\mathcal{A}(x_p,x_e)$.
	\end{prop}
	\begin{pf}
		In the case with zero capture radius and no obstacles,  $\mathcal{D}(x_p,x_e)=\mathcal{A}(x_p,x_e)$. We obtain the conclusion by \thref{shortestline_in_domainance}.\qed
	\end{pf}
	\begin{prop}\thlabel{uniq_C-oval}
		For any $x\in\overline{\mathcal{C}(x_p,x_e)}$, select a point $x'$ from the line segment whose endpoints are $x,x_e$. Then $x'\in\mathcal{C}(x_p,x_e)$.
	\end{prop}
	\begin{pf}
		Assume that there exists a point $x'$ on the line segment whose  endpoints are $x,x_e$ such that $x'\notin\mathcal{C}(x_p,x_e)$, i.e.$\left\| x'\right\|+ \left\| x_p\right\|-\alpha\left\| x'-xe\right\|\le 0 $. Then,
		\begin{align*}
			&\left\| x \right\|+ \left\| x_p\right\| \\
			\le &\left\| x-x'\right\| +\left\| x'\right\| +\left\| x_p\right\| \\
			<&\alpha\left\| x-x'\right\| +\left\| x'\right\| +\left\| x_p\right\| \\
			\le &\alpha\left\| x-x'\right\|+\alpha \left\| x'-x_e\right\| \\
			=&\alpha\left\| x-x_e\right\|.
		\end{align*}
		This contradicts $x\in\overline{\mathcal{C}(x_p,x_e)}$. \qed
	\end{pf}
	\begin{figure}[h]
		\centering
		\includegraphics[height=6cm]{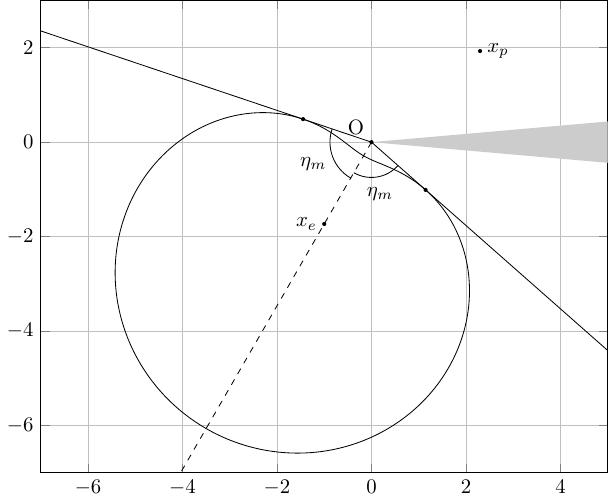}
		\caption{Tangent lines of $\partial\mathcal{C}(x_p,x_e)$ that pass through the origin.}
		\label{ztest0409_3} 
	\end{figure}
	$X$ in Section \ref{a_simple_case} is symmetric about the x-axis. Without loss of generality, assume that $x_p$ is in the upper half plane. Given $x_p,x_e$ with polar coordinate $(\rho_p,\theta_p),(\rho_e,\theta_e)$, take the two tangents of $\partial\mathcal{C}(x_p,x_e)$ through the origin. Let $\eta_m(x_p,x_e)$ denote the angles between the tangents and the line through $x_e$ and the origin, as shown in Fig.~\ref{ztest0409_3}. $\eta_m(x_p,x_e)$ can be abbreviated as $\eta_m$. According to \cite[Theorem 1]{garcia2021cooperative}, we obtain
	\begin{align*}
		\eta_m=\arccos\left( \frac{\sqrt{(\alpha^2-1)(\alpha^2\left\| x_e\right\|^2-\left\| x_p\right\|^2  )}-\left\| x_p\right\|}{\alpha^2 \left\| x_e\right\| }\right) 
	\end{align*}
	and 
	\begin{align*}
		\overline{\mathcal{C}(x_p,x_e)} \subseteq \mathcal{G}
		([\theta_e-\eta_m,\theta_e+\eta_m]) .
	\end{align*}
	We define a function on $\mathcal{G}
	([\theta_e-\eta_m,\theta_e+\eta_m]) $:
	\begin{equation}\label{eq_def_F}
		\begin{split}
			&f(x;x_p,x_e)\\
			=&
			\begin{cases}
				\left\| x-x_p\right\| -\alpha\left\| x-x_e\right\|  &x\in\mathcal{G}(I_1)\\
				\left\| x\right\| +\left\| x_p\right\| -\alpha\left\| x-x_e\right\| &x\in\mathcal{G}(I_2)
			\end{cases},
		\end{split}
	\end{equation}
	where
	\begin{equation*}
		\begin{split}
			&I_1=(-\infty,\theta_p+\pi]\cap[\theta_e-\eta_m,\theta_e+\eta_m],\\
			&I_2=(\theta_p+\pi,+\infty)\cap[\theta_e-\eta_m,\theta_e+\eta_m].
		\end{split}
	\end{equation*}
	Let
	\begin{equation}\label{set_def_F}
		\begin{split}
			&\mathcal{F}(x_p,x_e)\triangleq\\
			&\left\lbrace x\in\mathcal{G}([\theta_e-\eta_m,\theta_e+\eta_m]):f(x;x_p,x_e)>0\right\rbrace.
		\end{split}
	\end{equation}
	It is easy to obtain that
	\begin{equation}\label{set_def_partialF}
		\begin{split}
			&\partial\mathcal{F}(x_p,x_e)=\\
			&\left\lbrace x\in\mathcal{G}
			([\theta_e-\eta_m,\theta_e+\eta_m]) :f(x;x_p,x_e)=0\right\rbrace.
		\end{split}
	\end{equation}
	\begin{prop}\thlabel{D_in_F}
		$\partial\mathcal{D}(x_p,x_e)\subseteq\partial\mathcal{F}(x_p,x_e)$.
	\end{prop}
	\begin{pf}
		Synthesizing \eqref{eq_def_F} and \eqref{distance_of_shortest_path}, we obtain
		\begin{equation}\label{f_eq_d_L-ad_L}
			\begin{split}
				f(x;x_p,x_e)=d_L(x,x_p)-\alpha d_L(x,x_e),\\
				\forall x\in X\cap\mathcal{G}([\theta_e-\eta_m,\theta_e+\eta_m]).
			\end{split}
		\end{equation}
		It is easy to obtain that 
		\begin{align*}
			\partial\mathcal{A}(x_p,x_e)\subseteq \overline{\mathcal{C}(x_p,x_e)} \subseteq \mathcal{G}([\theta_e-\eta_m,\theta_e+\eta_m]), \\
			\partial\mathcal{C}(x_p,x_e)\subseteq \overline{\mathcal{C}(x_p,x_e)} \subseteq \mathcal{G}([\theta_e-\eta_m,\theta_e+\eta_m]).
		\end{align*}
		Synthesizing the above equations and \eqref{boundary=A+D}, we obtain
		\begin{align*}
			\partial\mathcal{D}(x_p,x_e)\subseteq \mathcal{G}
			([\theta_e-\eta_m,\theta_e+\eta_m]). 
		\end{align*}
		According to the definition of $\mathcal{D}(x_p,x_e)$,
		\begin{align*}
			\partial\mathcal{D}(x_p,x_e)\subseteq X.
		\end{align*}
		Thus, for any $x\in\partial\mathcal{D}(x_p,x_e)$,
		\begin{align*}
			f(x;x_p,x_e)=d_L(x,x_p)-\alpha d_L(x,x_e)=0,
		\end{align*}
		i.e. $x\in\partial\mathcal{F}(x_p,x_e)$. \qed
	\end{pf}
	\begin{prop}\thlabel{F-exist-uniq}
		For any $u_e \in \partial B_2(\boldsymbol{0},1)$, take a ray from the point $x_e$ along the vector $u_e$. Then, there exists a unique intersection point of the ray and $\mathcal{F}(x_p,x_e)$.
	\end{prop}
	\begin{pf}
		We first prove the existence. It is easy to obtain that $\lim\limits_{\left\| x\right\| \rightarrow +\infty}(\left\| x\right\|+\left\| x_p\right\| -\alpha\left\| x-x_e\right\|)=-\infty $. Thus, $\partial\mathcal{C}(x_p,x_e)$ is bounded. Pick $u_e \in \partial B_2(\boldsymbol{0},1)$ and take a ray from the point $x_e$ along the vector $u_e$. There exists an intersection point of the ray and $\partial\mathcal{C}(x_p,x_e)$, which is denoted by $x_0$. If $f(x_0;x_p,x_e)=0$, $x_0$ is the intersection point of the ray and $\partial\mathcal{F}(x_p,x_e)$. If $f(x_0;x_p,x_e)\neq 0$,
		\begin{equation*}
			\begin{split}
				f(x_0;x_p,x_e)&=\left\| x_0-x_p\right\| -\alpha\left\| x_0-x_e\right\| \\
				&< \left\| x_0\right\| +\left\| x_p\right\| -\alpha\left\| x_0-x_e\right\|=0.
			\end{split}
		\end{equation*}
		We also have $f(x_e;x_p,x_e)>0$. By the intermediate value theorem, there exists a point $x_0'$ on the line segment whose endpoints are $x_0$ and $x_e$ such that $f(x_0';x_p,x_e)=0$. $x_0'$ is the intersection point of the ray and $\partial\mathcal{F}(x_p,x_e)$. Existence has been proven.\\
		Next, we prove the uniqueness. For the sake of contradiction, assume that there are two different intersection points of the ray and $\partial\mathcal{F}(x_p,x_e)$, which are denoted as $x_1,x_2$. Without loss of generality, let $\left\| x_1-x_e\right\| >\left\|x_2-x_e\right\|$. According to \thref{uniq_A-circle} and \thref{uniq_C-oval}, one of $x_1,x_2$ is on $\partial\mathcal{A}(x_p,x_e)$, the other is on $\partial\mathcal{C}(x_p,x_e)$. \\
		If $x_1\in\partial\mathcal{A}(x_p,x_e)$ and $x_2\in\partial\mathcal{C}(x_p,x_e)$, we obtain
		\begin{align*}
			\left\| x_2\right\| +\left\| x_p\right\| -\alpha\left\| x_2-x_e\right\| =0.
		\end{align*}
		According to \thref{uniq_A-circle},
		\begin{align*}
			\left\| x_2- x_p\right\| -\alpha\left\| x_2-x_e\right\| >0.
		\end{align*}
		Then,
		\begin{align*}
			\left\| x_2- x_p\right\|>\left\| x_2\right\| +\left\| x_p\right\|,
		\end{align*}
		which contradicts the triangle inequality of norms.\\ 
		If $x_1\in\partial\mathcal{C}(x_p,x_e)$ and $x_2\in\partial\mathcal{A}(x_p,x_e)$, 
		due to \thref{uniq_A-circle} and $\lim\limits_{\left\| x\right\| \rightarrow +\infty }(\left\| x- x_p\right\| -\alpha\left\| x_2-x_e\right\| )=-\infty$, 
		\begin{align*}
			&\left\| x- x_p\right\| -\alpha\left\| x-x_e\right\|<0,\\
			&\forall x\in\left\lbrace \lambda x_1+(1-\lambda)x_2:\lambda\in(0,1)\right\rbrace .
		\end{align*}
		According to \thref{uniq_C-oval},
		\begin{align*}
			&\left\| x\right\| +\left\|  x_p\right\| -\alpha\left\| x-x_e\right\|>0,\\
			&\forall x\in\left\lbrace \lambda x_1+(1-\lambda)x_2:\lambda\in(0,1)\right\rbrace. 
		\end{align*}
		Thus,
		\begin{equation}\label{app_3}
			\begin{split}
				&\left\| x\right\| +\left\|  x_p\right\| >\left\| x- x_p\right\|,\\
				&\forall x\in\left\lbrace \lambda x_1+(1-\lambda)x_2:\lambda\in(0,1)\right\rbrace. 
			\end{split}
		\end{equation}
		According to the definition of $\mathcal{F}(x_p,x_e)$ \eqref{set_def_F}, we obtain $x_1 \in \mathcal{G}(I_2)$ and $x_2\in \mathcal{G}(I_1)$. Thus, the line segment connecting $x_1$ and $x_2$ intersects with the boundary of $\mathcal{G}(I_1)$ and $\mathcal{G}(I_2)$:
		\begin{align*}
			\left\lbrace (\rho\cos\theta,\rho\sin\theta):\rho\ge0,\theta=\theta_e+\pi \right\rbrace. 
		\end{align*}
		Let $x_0$ denote the intersection point. $x_0$ satisfies $\left\| x_0- x_p\right\|=\left\| x_0\right\| +\left\|  x_p\right\|$. Due to \eqref{app_3}, $x_0 \notin \left\lbrace \lambda x_1+(1-\lambda)x_2:\lambda\in(0,1)\right\rbrace$. According to \thref{uniq_A-circle} and \thref{uniq_C-oval}, $x_0$ cannot be $x_1$ or $x_2$, and a contradiction arises. 
		\qed
		\begin{rem}\thlabel{F_closed_curve}
			\thref{F-exist-uniq} established a bijection between $\partial\mathcal{F}(x_p,x_e)$ and $\partial B_2(\boldsymbol{0},1)$. The map is continuously differentiable because $f$ is continuously differentiable. This indicates that $\partial\mathcal{F}(x_p,x_e)$ is a continuously differentiable, simple, closed curve.
		\end{rem}
	\end{pf}
	\section{Proof of \thref{lemma_of_increament}}\label{proof_lemma_of_increament}
	\begin{pf}
		We will prove this proposition in 4 cases.\\
		Case 1: $x_p,x_e$ are visible to each other and $x\in \partial \mathcal{A}(x_p,x_e)$. We have 
		\begin{align*}
			-\partial_{1}d_L(x_p,x_e)=\frac{x_e-x_p}{\left\| x_e-x_p \right\| },\\
			- \partial_{2}d_L(x,x_p)=\frac{x-x_p}{\left\| x-x_p \right\| }.
		\end{align*}
		\begin{figure}[h]
			\centering
			\includegraphics[height=6cm]{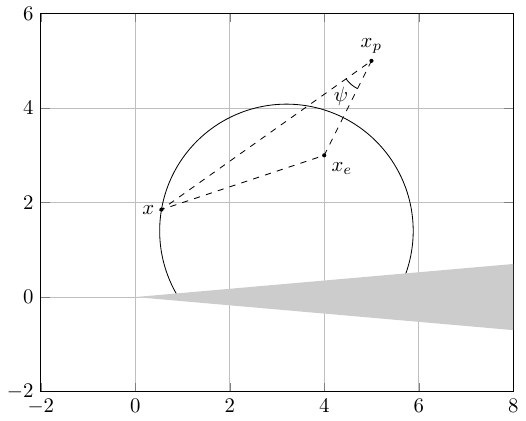}
			\caption{In the proof of \thref{lemma_of_increament}, $x_p,x_e$ are visible to each other and $x \in \partial \mathcal{A}(x_p,x_e)$.}
			\label{zincrement001} 
		\end{figure}\\
		As shown in Fig.~\ref{zincrement001}, let $\psi$ denote the angle between vectors $x_e-x_p$ and $x-x_p$. By \cite[Theorem 1]{garcia2021cooperative},
		\begin{align*}
			\psi\le \arcsin\frac{1}{\alpha} <\frac{\pi}{2}.
		\end{align*}
		Thus,
		\begin{align*}
			\partial_{1}d_L(x_p,x_e) \cdot  \partial_{2}d_L(x,x_p)=\cos\psi>0.
		\end{align*}
		Case 2: $x_p,x_e$ are visible to each other and $x\in \partial \mathcal{C}(x_p,x_e)$. We have
		\begin{align*}
			-\partial_{1}d_L(x_p,x_e)=\frac{x_e-x_p}{\left\| x_e-x_p \right\| },\\
			- \partial_{2}d_L(x,x_p)=-\frac{x_p}{\left\| x_p \right\| }.
		\end{align*}
		\begin{figure}[h]
			\centering
			\includegraphics[height=6cm]{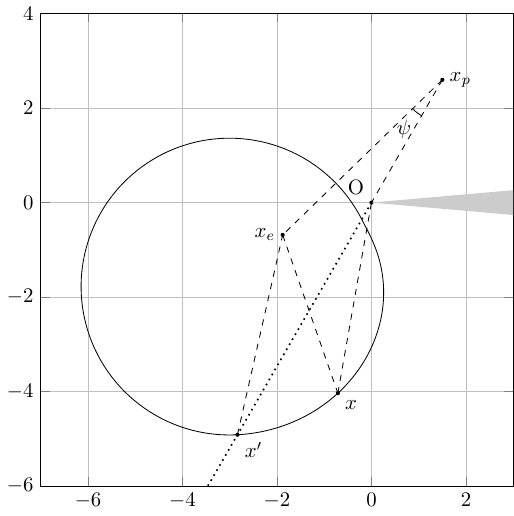}
			\caption{In the proof of \thref{lemma_of_increament}, $x_p,x_e$ are visible to each other and $x\in \partial \mathcal{C}(x_p,x_e)$.}
			\label{zincrement002}
		\end{figure}\\
		As shown in Fig.~\ref{zincrement002}, there exist two intersection points of $\partial \mathcal{A}(x_p,x_e)$ and the line connecting $x_p$ and the origin. Let $x'$ denote the farther intersection point from $x_p$. We have
		\begin{align*}
			- \partial_{2}d_L(x',x_p)=\frac{x'-x_p}{\left\| x'-x_p \right\| }=-\frac{x_p}{\left\| x_p \right\| }=- \partial_{2}d_L(x,x_p).
		\end{align*}
		Let $\psi$ denote the angle between vectors $x_e-x_p$ and $x'-x_p$. Similarly to Case 1, by \cite[Theorem 1]{garcia2021cooperative},
		\begin{align*}
			\psi\le \arcsin\frac{1}{\alpha} <\frac{\pi}{2}.
		\end{align*}
		Thus,
		\begin{align*}
			&\partial_{1}d_L(x_p,x_e) \cdot  \partial_{2}d_L(x,x_p),\\
			=&\partial_{1}d_L(x_p,x_e) \cdot  \partial_{2}d_L(x',x_p)=\cos\psi>0.
		\end{align*}
		Case 3: $x_p,x_e$ are not visible to each other and $x \in \partial \mathcal{A}(x_p,x_e)$. We have
		\begin{align*}
			-\partial_{1}d_L(x_p,x_e)=-\frac{x_p}{\left\| x_p \right\| },\\
			- \partial_{2}d_L(x,x_p)=\frac{x-x_p}{\left\| x-x_p \right\| }.
		\end{align*}
		\begin{figure}[h]
			\centering
			\includegraphics[height=6cm]{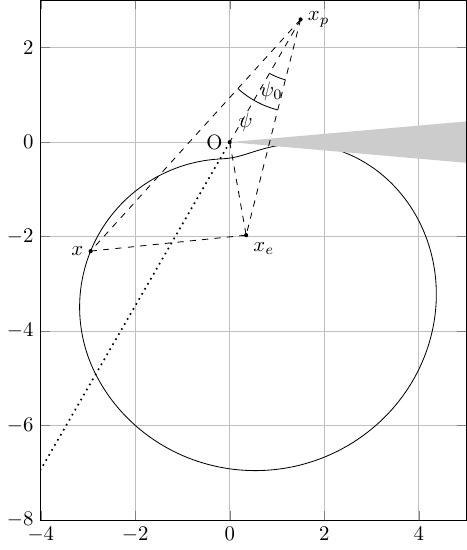}
			\caption{In the proof of \thref{lemma_of_increament}, $x_p,x_e$ are not visible to each other and $x \in \partial \mathcal{A}(x_p,x_e)$.}
			\label{zincrement003} 
		\end{figure}\\
		As shown in Fig.~\ref{zincrement003}, let $\psi$ denote the angle between vectors $x-x_p$ and $x_e-x_p$, and let $\psi_0$ denote the angle between vectors $-x_p$ and $x_e-x_p$. It is easy to obtain that the origin and $x$ are on the same side of the line connecting $x_p,x_e$. According to \cite[Theorem 1]{garcia2021cooperative}, we obtain
		\begin{align*}
			\psi_0\le\psi\le \arcsin\frac{1}{\alpha} <\frac{\pi}{2}.
		\end{align*}
		Thus,
		\begin{align*}
			&\partial_{1}d_L(x_p,x_e) \cdot  \partial_{2}d_L(x,x_p)\\
			=&\cos(\psi-\psi_0)\ge \cos\psi 
			>0.			
		\end{align*}
		Case 4: $x_p,x_e$ are not visible to each other and $x \in \partial \mathcal{C}(x_p,x_e)$. We have
		\begin{align*}
			\partial_{1}d_L(x_p,x_e)= \partial_{2}d_L(x,x_p)=\frac{x_p}{\left\| x_p \right\|}.
		\end{align*}
		Thus,
		\begin{align*}
			\partial_{1}d_L(x_p,x_e) \cdot  \partial_{2}d_L(x,x_p)=1>0.
		\end{align*}
		The proof is completed. \qed
	\end{pf}
\end{document}